\documentclass[reqno,12pt]{amsart}
\usepackage{amsmath,amsthm,amscd,amssymb}
\newcommand{\bbC}{{\mathbb{C}}}
\newcommand{\bbD}{{\mathbb{D}}}

\newcommand{\bbO}{{\mathbb{O}}}

\newcommand{\bbR}{{\mathbb{R}}}
\newcommand{\bbS}{{\mathbb{S}}}
\newcommand{\bbT}{{\mathbb{T}}}
\newcommand{\bbU}{{\mathbb{U}}}
\newcommand{\bbZ}{{\mathbb{Z}}}

\newcommand{\calB}{{\mathcal B}}
\newcommand{\calC}{{\mathcal C}}
\newcommand{\calE}{{\mathcal E}}
\newcommand{\calF}{{\mathcal F}}
\newcommand{\calG}{{\mathcal G}}
\newcommand{\calH}{{\mathcal H}}
\newcommand{\calI}{{\mathcal I}}
\newcommand{\calJ}{{\mathcal J}}

\newcommand{\calL}{{\mathcal L}}
\newcommand{\calM}{{\mathcal M}}
\newcommand{\calR}{{\mathcal R}}
\newcommand{\calS}{{\mathcal S}}
\newcommand{\calO}{{\mathcal O}}

\newcommand{\lb}{\label}
\newcommand{\f}{\frac}

\newcommand{\ol}{\overline}
\newcommand{\ti}{\tilde  }
\newcommand{\wti}{\widetilde  }
\newcommand{\llangle}{\langle\!\langle}
\newcommand{\rrangle}{\rangle\!\rangle}

\newcommand{\tr}{\text{\rm{Tr}}}

\newcommand{\bdone}{{\boldsymbol{1}}}

\newcommand{\ess}{\text{\rm{ess}}}
\newcommand{\ac}{\text{\rm{ac}}}

\newcommand{\capL}{\text{\rm{L}}}
\newcommand{\capR}{\text{\rm{R}}}

\newcommand{\supp}{\text{\rm{supp}}}
\newcommand{\Span}{\text{\rm{span}}}

\newcommand{\bi}{\bibitem}

\newcommand{\beq}{\begin{equation}}
\newcommand{\eeq}{\end{equation}}
\newcommand{\ba}{\begin{align}}
\newcommand{\ea}{\end{align}}
\newcommand{\veps}{\varepsilon}

\DeclareMathOperator{\Real}{Re}
\DeclareMathOperator{\ran}{Ran}
\DeclareMathOperator{\Ima}{Im}

\DeclareMathOperator*{\slim}{s-lim}

\let\det=\undefined\DeclareMathOperator{\det}{det}
\allowdisplaybreaks
\numberwithin{equation}{section}

\newtheorem{theorem}{Theorem}[section]

\newtheorem{proposition}[theorem]{Proposition}
\newtheorem{lemma}[theorem]{Lemma}
\newtheorem{corollary}[theorem]{Corollary}
\theoremstyle{definition}

\newtheorem{example}[theorem]{Example}

\theoremstyle{remark}
\newtheorem*{remark}{Remark}
\newtheorem*{remarks}{Remarks}

\newcommand{\abs}[1]{\lvert#1\rvert}

\newcounter{smalllist}
\newenvironment{SL}{\begin{list}{{\rm\roman{smalllist})}}{%
\setlength{\topsep}{0mm}\setlength{\parsep}{0mm}\setlength{\itemsep}{0mm}%
\setlength{\labelwidth}{2em}\setlength{\leftmargin}{2em}\usecounter{smalllist}%
}}{\end{list}}

\begin{document}

\title{CMV Matrices: Five Years After}

\author[B.~Simon]{Barry Simon*}

\thanks{$^*$ Mathematics 253-37, California Institute of Technology, Pasadena, CA 91125, USA.
E-mail: bsimon@caltech.edu. Supported in part by NSF grant DMS-0140592}
\thanks{Submitted to the Proceedings of the W.~D.~Evans' 65th Birthday Conference}

\date{March 1, 2006}

\begin{abstract} CMV matrices are the unitary analog of Jacobi matrices;
we review their general theory.
\end{abstract}

\maketitle

\section{Introduction} \lb{s1}

\footnotesize
\begin{flushleft}
{\sffamily {\it The Arnold Principle:} If a notion bears a personal
name, then this name is not the name of the inventor. \\
{\it The Berry Principle:} The Arnold Principle is applicable
to itself. }

\hfill{\slshape --- V.~I.~Arnold} \\
\hfill{\slshape On Teaching Mathematics, 1997 \cite{Arnold}} \\
\hfill{(Arnold says that Berry formulated these principles.)}

\end{flushleft}

\normalsize

\bigskip

In 1848, Jacobi \cite{Jaco} initiated the study of quadratic forms
$J(x_1, \dots, x_n) = \sum_{k=1}^n b_k x_k^2 + 2 \sum_{k=1}^{n-1}
a_k x_k x_{k+1}$, that is, essentially $n\times n$ matrices of the form
\begin{equation} \lb{1.1}
J= \begin{pmatrix}
b_1 & a_1 & 0   & \dots & 0 \\
a_1 & b_2 & a_2 & \dots & 0 \\
0   & a_2 & b_3 & \dots & 0 \\
\vdots & \vdots  & \vdots & \ddots & \vdots \\
0 &\dots & \dots & a_{n-1} &  b_n
\end{pmatrix}
\end{equation}
and found that the eigenvalues of $J$ were the zeros of the denominator
of the continued fraction
\begin{equation} \lb{1.2}
\cfrac{1}{b_1 - z - \cfrac{a_1^2}{b_2 -z - \cfrac{a_2^2}{\cdots} }}
\end{equation}
In the era of the birth of the spectral theorem, Toeplitz \cite{T1910},
Hellinger--Toeplitz \cite{HellTop}, and especially Stone \cite{Stone} realized
that Jacobi matrices were universal models of selfadjoint operators, $A$, with
a cyclic vector, $\varphi_0$.

To avoid technicalities, consider the case where $A$ is bounded, and suppose
initially that $\calH$ is infinite-dimensional. By cyclicity, $\{A^k\varphi_0\}_{k=0}^\infty$
are linearly independent, so by applying Gram--Schmidt to $\varphi_0, A\varphi_0,
A^2 \varphi_0,\dots$, we get polynomials $p_j(A)$ of degree exactly $j$ with positive
leading coefficients so that
\begin{equation} \lb{1.3}
\varphi_j = p_j(A) \varphi_0
\end{equation}
are an orthonormal basis for $\calH$. By construction,
\[
\varphi_j \perp \varphi_0, A\varphi_0, \dots, A^{j-1} \varphi_0
\]
so
\begin{equation} \lb{1.4}
\langle\varphi_j, A\varphi_k\rangle =0 \qquad j\geq k+2
\end{equation}
Because $A$ is selfadjoint, we see $\langle \varphi_j, A\varphi_k\rangle =0$
also if $j\leq k-2$. Thus, the matrix $\langle\varphi_j, A\varphi_k\rangle$
has exactly the form \eqref{1.1} where $a_j >0$ (since $p_j(A)$ has leading
positive coefficient).

Put differently, for all $A,\varphi_0$, there is a unitary $U:\calH\to\ell^2$
(given by Fourier components in the $\varphi_j$ basis), so $U\!AU^{-1}$ has
the form $J$ and $\varphi_0 = (1, 0, 0, \dots)^t$. The Jacobi parameters,
$\{a_n,b_n\}_{n=1}^\infty$, are intrinsic, which shows there is exactly one $J$
(with $\varphi_0 = (1, 0, 0, \dots)^t$) in the unitary equivalence class
of $(A,\varphi_0)$.

There is, of course, another way of describing unitary invariants for $(A, \varphi_0)$:
the spectral measure $d\mu$ defined by
\begin{equation} \lb{1.5}
\int x^n \, d\mu(x) = \langle\varphi_0, A^n\varphi_0\rangle
\end{equation}
There is a direct link from $d\mu$ to the Jacobi parameters: the $p_j(x)$ are
orthonormal polynomials associated to $d\mu$, and the Jacobi parameters are
associated to the three-term recursion relation obeyed by the $p$'s:
\begin{equation} \lb{1.6}
xp_j(x) = a_{j+1} p_{j+1} + b_{j+1} p_j(x) + a_j p_{j-1}(x)
\end{equation}
(where $p_{-1}\equiv 0$).

Here we are interested in the analog of these structures for unitary matrices.
We begin by remarking that for a general normal operator, $N$, the right form of
cyclicity is that $\{N^k (N^*)^\ell\varphi_0\}_{k,\ell=0}^\infty$ is total.
Since $A=A^*$, only $\{A^k \varphi_0\}_{k=0}^\infty$ enters. Since $U^*=U^{-1}$,
for unitaries $U^k (U^*)^\ell=U^{k-\ell}$ and the right notion of cyclicity
is that $\{U^k\varphi_0\}_{k=-\infty}^\infty$ is total.

Some parts of the above four-fold equivalence:
\begin{SL}
\item[(1)] unitary equivalence classes of $(A,\varphi_0)$
\item[(2)] spectral measures, that is, probability measures $d\mu$ on $\bbR$
with bounded support and infinite support
\item[(3)] Jacobi parameters
\item[(4)] Jacobi matrices
\end{SL}
are immediate for the unitary case. Namely, (1) $\Leftrightarrow$ (2) holds since
there is a spectral theorem for unitaries, and so, a one-one correspondence between unitary
equivalence classes of $(U,\varphi_0)$ on infinite-dimensional spaces and
probability measures on $\partial\bbD$ ($\bbD=\{z\mid\abs{z}<1\}$) with infinite
support.

More subtle is the analog of Jacobi parameters. Starting from such a probability
measure on $\partial\bbD$, one can form the monic orthogonal polynomials $\Phi_n(z)$
and find (see Szeg\H{o} \cite{Szbk}; see also Section~1.5 of \cite{OPUC1})
$\{\alpha_n\}_{n=0}^\infty\in \bbD^\infty$, so
\begin{equation} \lb{1.7}
z\Phi_n(z) =\Phi_{n+1}(z) + \bar\alpha_n z^n\, \ol{\Phi_n (1/\bar z)}
\end{equation}
While Verblunsky \cite{V35} defined the $\alpha_n$ in a different (albeit
equivalent) way, he proved a theorem (called Verblusnky's theorem in \cite{OPUC1};
see also \cite{1Foot}) that says this map $d\mu\to\{\alpha_n\}_{n=0}^\infty$ is
one-one and onto all of $\bbD^\infty$, so (1)--(3) for the unitary case have been
well understood for sixty-five years.

Surprisingly, (4) (i.e., the canonical matrix form for unitaries) is of much
more recent vintage. The key paper by Cantero, Moral, and Vel\'azquez \cite{CMV}
was submitted in April 2001 --- so we are witnessing five years of study in
the area --- it is reviewing these developments that is the main scope of this
review article. Spectral theory of differential and difference operators has
been an enduring theme of Des Evans' research and I am pleased to dedicate
this review to him.

There is an ``obvious" matrix to try, namely, $\calG_{k\ell}=
\langle\varphi_k, z\varphi_\ell\rangle$ with $\varphi_k = \Phi_k/\|\Phi_k\|$
the orthonormal polynomials. This GGT matrix (as it is named in \cite{OPUC1};
see Section~\ref{s10} below) has two defects. First, $\{\varphi_k\}_{k=0}^\infty$
is a basis if and only if $\sum_{n=0}^\infty \abs{\alpha_n}^2=\infty$, and if
it is not, $\calG_{k\ell}$ is not unitary and is not conjugate to multiplication
by $z$ (in that case, one can look at the minimal dilation of $\calG$, which
is discussed in \cite{Con84,OPUC1}). Second, it obeys \eqref{1.4}, but in general,
$\langle\varphi_j, U^*\varphi_k\rangle\neq 0$ for all $j\geq k+1$, that is,
$\calG$ is not of finite width measured from the diagonal. CMV \cite{CMV}
has the following critical ideas:
\begin{SL}
\item[(a)] The basis $\chi_k$ obtained by orthonormalizing $1,z,z^{-1}, z^2,
z^{-2}, \dots$ can be written in terms of $\varphi_\ell(z)$,
$\ol{\varphi_\ell (1/\bar z)}$, and powers of $z$.

\item[(b)] The matrix $\calC_{k\ell}= \langle \chi_k, z\chi_\ell\rangle$
is unitary and five-diagonal.

\item[(c)] $\calC$ can be factorized into $\calC=\calL\calM$ where $\calL$ is
a direct sum of $2\times 2$ unitary matrices and $\calM$ the direct sum of a
single $1\times 1$ and $2\times 2$ matrices.
\end{SL}

It turns out that these key ideas appeared about ten years earlier in the
numeric matrix literature (still, of course, much later than the 1930's
resolution of (1)--(3)). Intimately related to this history is what we
will call the AGR factorization in Section~\ref{s11} --- the ability to write $\calG$
in the case of $n\times n$ matrices as a product $\wti\Theta_0\dots\wti\Theta_{n-1}
\wti{\wti\Theta}_{n-1}$ of matrices with a single $2\times 2$ block placed in
$\bdone$ and a finite matrix which is diagonal, differing from $\bdone$ in a single
place (see Section~\ref{s11} for details).

In 1986, Ammar, Gragg, and Reichel \cite{Athens} found the AGR factorization for
orthogonal matrices --- here the $\alpha_j$ are real and the $\Theta (\alpha_j)$ are
reflections, so the AGR factorization can be viewed as an iteration of a Householder
algorithm. In this paper, they also had a proof of the $\calL\calM$ factorization for
this case. This proof (a variant of which appears in Section~\ref{s10}), which works
in general to go from the AGR factorization of the GGT matrix to an $\calL\calM$
factorization, was only given in the orthogonal case since they did not yet have
the AGR factorization for general unitaries.

In 1988 (published 1991), AGR \cite{AGR} extended the AGR factorization to the general
unitary case and realized the connection to Szeg\H{o} recursion. While they could
have proven an $\calL\calM$ factorization from this using the method in \cite{Athens},
they did not and the general $\calL\calM$ factorization only appeared in \cite{BGE}.

In 1991, Bunse-Gerstner and Elsner \cite{BGE} found the $\calL\calM$ factorization for
a general finite unitary and noted it was a five-diagonal representation. Watkins
\cite{Watk} codified and cleaned up those results and emphasized the connection to OPUC
and found a proof of Szeg\H{o} recursion from the $\calL\calM$ factorization. Virtually
all the main results from \cite{CMV} are already in Watkins \cite{Watk}.

We will continue to use the name CMV matrices, in part because the analytic revolution
we discuss here was ushered in by their work and in part because the name
has been used now in many, many publications.

Here is a summary of the rest of this review. Section~\ref{s2} presents the basics,
essentially notation and (a)--(c) above. Section~\ref{s3} discusses ``other" CMV
matrices. In particular, we consider two kinds of finite variants. In the selfadjoint
case, restricting the matrix by taking the first $n$ rows and columns preserves
selfadjointness but the analog for unitaries does not, and we have both the nonunitary
cutoff CMV matrices obtained from the first $n$ rows and columns and the unitary finite
 CMV matrices which are models of finite unitary matrices with a cyclic vector.
Section~\ref{s4} discusses CMV matrices for matrix-valued measures --- something
that is new here. Section~\ref{s5} discusses the effect on Verblunsky coefficients
of rank one multiplication perturbations, and Section~\ref{s6} the formula for
the resolvent of the CMV matrices, the analog of well-known Green's function
formulae for Jacobi matrices. Sections~\ref{s7} and \ref{s9} discuss perturbation
results, and Section~\ref{s8} a general theorem on the essential spectrum of
CMV matrices. Section~\ref{s10} discusses the AGR factorization discussed
above as preparation for the Killip--Nenciu discussion of five-diagonal models
for $\beta$-distribution of eigenvalues, the subject of Section~\ref{s11}.
Section~\ref{s12} discusses the defocusing AL flows, which bear the same
relation to CMV matrices as Toda flows do to Jacobi matrices. Finally,
Section~\ref{s13} discusses a natural reduction of CMV matrices to a
direct sum of two Jacobi matrices when all Verblunsky coefficients are real.

We do not discuss the use of CMV matrices to compute the zeros of OPUC. These
zeros are the eigenvalues of the cutoff CMV matrix. We note that this method
of computing zeros was used in the recent paper of Mart\'inez-Finkelshtein,
McLaughlin, and Saff \cite{MFMS}. Numerical aspects of CMV matrices deserve
further study.

While this is primarily a review article, there are numerous new results, including:
\begin{SL}
\item[(1)] an analysis of what matrices occur as cutoff CMV matrices (Section~\ref{s3})
\item[(2)] an analysis following Watkins \cite{Watk} of the $\calL\calM$ factorization
without recourse to Szeg\H{o} recursion (Section~\ref{s3})
\item[(3)] the basics of CMV matrices for matrix-valued measures (Section~\ref{s4})
\item[(4)] a new proof of AGR factorization using intermediate bases (Section~\ref{s10})
\item[(5)] a new trace class estimate for GGT matrices that relies on AGR factorization
(Section~\ref{s10})
\item[(6)] a reworked proof of the Killip--Nenciu \cite{KN} theorem on the measure that
Haar measure on $\bbU(n)$ induces on Verblunsky coefficients (Section~\ref{s11})
\item[(7)] an argument of AGR is made explicit and streamlined to go from AGR to
$\calL\calM$ factorization (Section~\ref{s10})
\end{SL}

\medskip
It is a pleasure to thank M.~Cantero, P.~Deift, L.~Golinskii, F.~Marcell\'an,
A.~Mart\'inez-Finkelshtein, L.~Moral, I.~Nenciu, P.~Nevai, L.~Vel\'azquez, and
D.~Watkins for useful input.

\section{CMV Matrices: The Basics} \lb{s2}

In this section, we define the CMV basis, the CMV matrix,
and the $\calL\calM$ factorization.

CMV matrices can be thought of in terms of unitary matrices or OPUC. We start
with the OPUC point of view. A measure $d\mu$ in $\partial\bbD$ is called nontrivial
if it is not supported on a finite set; equivalently, if every polynomial,
which is not identically zero, is nonzero in $L^2 (\partial\bbD, d\mu)$.
Then one can define {\it orthonormal polynomials}, $\varphi_n(z)$ (or $\varphi_n
(z,d\mu))$, by
\begin{alignat}{2}
&\text{(i)} \qquad && \varphi_n(z) = \kappa_n z^n + \text{lower order}; \quad
\kappa_n >0 \lb{2.1} \\
&\text{(ii)} \qquad && \varphi_n \perp \{1,z,z^2, \dots, z^{n-1}\} \lb{2.2}
\end{alignat}
We define the monic polynomials $\Phi_n(z)$ by $\Phi_n(z) =\varphi_n(z)/\kappa_n$.

The {\it Szeg\H{o} dual} is defined by
\begin{equation} \lb{2.3}
P_n^*(z) = z^n \, \ol{P_n (1/\bar z)}
\end{equation}
that is,
\begin{equation} \lb{2.4}
P_n(z) = \sum_{j=0}^n c_j z^j \Rightarrow P_n^*(z) = \sum_{j=0}^n
\bar c_{n-j} z^j
\end{equation}
The symbol $^*$ is $n$-dependent and is sometimes applied to polynomials of degree at
most $n$, making the notation ambiguous!

Then there are constants $\{\alpha_n\}_{n=0}^\infty$ in $\bbD$, called {\it Verblunsky
coefficients}, (sometimes we will write $\alpha_n (d\mu)$) so that
\begin{equation} \lb{2.5}
\rho_n\,\, \varphi_{n+1}(z) = z\varphi_n(z) -\bar\alpha_n \varphi_n^*(z)
\end{equation}
where
\begin{equation} \lb{2.6x}
\rho_n = (1-\abs{\alpha_n}^2)^{1/2}
\end{equation}
Moreover, $\mu\to \{\alpha_n\}_{n=0}^\infty$ sets up a one-one correspondence between
nontrivial measures on $\partial\bbD$ and points of $\bbD^\infty$ (as we will show below).
\eqref{2.5} (called Szeg\H{o} recursion after \cite{Szbk}) and this one-one
correspondence are discussed in \cite{OPUC1,OPUC2}; see also \cite{1Foot}.
Applying ${}^*$ for $P_{n+1}$ to \eqref{2.5}, we get
\begin{equation} \lb{2.6}
\rho_n\,\, \varphi_{n+1}^* (z) = \varphi_n^*(z) -\alpha_n z \varphi_n(z)
\end{equation}

If one defines (of course, $\rho = (1-\abs{\alpha}^2)^{1/2}$)
\begin{equation} \lb{2.7}
A(\alpha) = \f{1}{\rho} \begin{pmatrix} z & -\bar\alpha \\
-\alpha z & 1 \end{pmatrix}
\end{equation}
then \eqref{2.5}/\eqref{2.6} can be written
\begin{equation} \lb{2.8}
\binom{\varphi_{n+1}}{\varphi_{n+1}^*} = A(\alpha_n) \binom{\varphi_n}{\varphi_n^*}
\end{equation}

Since $\det(A) =z$, we have
\begin{equation} \lb{2.9}
A(\alpha)^{-1} = \f{1}{\rho z}
\begin{pmatrix}
1 & \bar\alpha \\
\alpha z & z
\end{pmatrix}
\end{equation}
and thus
\begin{align}
\rho_n\,\, \varphi_n(z) &= \f{\varphi_{n+1}(z) + \bar\alpha_n \varphi_{n+1}^*(z)}{z} \lb{2.10} \\
\rho_n\,\, \varphi_n^*(z) &= (\varphi_{n+1}^*(z) + \alpha_n \varphi_{n+1}(z)) \lb{2.11}
\end{align}

Introduce the notation $[y_1, \dots, y_k]$ for the span of the vectors $y_1, \dots, y_k$
and $P_{[y_1, \dots, y_k]}$ for the projection onto the space $[y_1, \dots, y_k]$.
For $x\notin [y_1, \dots, y_k]$, define
\begin{equation} \lb{2.12}
[x;y_1, \dots, y_k] = \f{(1-P_{[y_1, \dots, y_k]})x}{\|(1-P_{[y_1, \dots, y_k]})x\|}
\end{equation}
the normalized projection of $x$ onto $[y_1, \dots, y_k]^\perp$, that is, the result of
adding $x$ to a Gram--Schmidt procedure.

We define $\pi_n$ to be $P_{[1,\dots, z^{n-1}]}$.

By the definition of $\varphi_n$ and the fact that ${}^*$ is anti-unitary on $\ran \pi_n$ and
takes $z^j$ to $z^{n-j}$, we have
\begin{equation} \lb{2.13}
\varphi_n = [z^n; 1, \dots, z^{n-1}] \qquad
\varphi_n^* = [1;z, \dots, z^n]
\end{equation}

With this notation out of the way, we can define the {\it CMV basis\/} $\{\chi_n\}_{n=0}^\infty$
and {\it alternate CMV basis\/} $\{x_n\}_{n=0}^\infty$ as the Laurent polynomials (i.e.,
polynomials in $z$ and $z^{-1}$) obtained by applying Gram--Schmidt to $1, z, z^{-1}, z^2,
z^{-2}, \dots$ and $1,z^{-1}, z,z^{-2}, z^2, \dots$, that is, for $k=0,1,\dots$,
\begin{alignat}{2}
\chi_{2k} &= [z^{-k}; 1,z,\dots, z^{-k+1}, z^k]
& \qquad \chi_{2k-1} &= [z^k; 1,z,\dots, z^{k-1}, z^{-k+1}] \lb{2.14} \\
x_{2k} &= [z^k; 1,z^{-1}, \dots, z^{k-1}, z^{-k}]
& \quad x_{2k-1} &= [z^{-k}; 1,z^{-1}, \dots, z^{-k+1}, z^{k-1}] \lb{2.15}
\end{alignat}
So, in particular, as functions in $L^2 (\partial\bbD, d\mu)$,
\begin{equation} \lb{2.16}
x_n = \bar\chi_n
\end{equation}
and as Laurent polynomials,
\begin{equation} \lb{2.17}
x_n(z) = \ol{\chi_n (1/\bar z)}
\end{equation}

As realized by CMV \cite{CMV}, the $\{\chi_n\}_{n=0}^\infty$ and $\{x_n\}_{n=0}^\infty$
are always a basis of $L^2 (\partial\bbD, d\mu)$ since the Laurent polynomials are dense
on $C(\partial\bbD)$, while $\{\varphi_n\}_{n=0}^\infty$ may or may not be a basis (it
is known that this is a basis if and only if $\sum_n \abs{\alpha_n}^2 =\infty$;
see Theorem~1.5.7 of \cite{OPUC1}). On the other hand, the $\chi$ and $x$ bases can be
expressed in terms of $\varphi$ and $\varphi^*$ by \eqref{2.13} and the fact that
multiplication by $z^\ell$ is unitary. For example,
\begin{align*}
x_{2k} &= z^{-k} [z^{2k}; z^k, z^{k-1}, \dots, z^{2k-1}, 1] \\
&= z^{-k} [z^{2k}; 1, \dots, z^{2k-1}] \\
&= z^{-k} \varphi_{2k}(z)
\end{align*}

The full set is
\begin{alignat}{2}
\chi_{2k}(z) &= z^{-k} \varphi_{2k}^*(z) & \qquad
\chi_{2k-1}(z) &= z^{-k+1} \varphi_{2k-1}(z) \lb{2.18} \\
x_{2k}(z) &= z^{-k}\varphi_{2k}(z) & \qquad
x_{2k-1}(z) & = z^{-k} \varphi_{2k-1}^*(z) \lb{2.19}
\end{alignat}

Since $\chi$ and $x$ are bases, the matrices of multiplication by $z$ in
these bases are unitary. So we have the unitary matrices
\begin{align}
\calC_{m\ell} &= \langle \chi_m, z\chi_\ell\rangle \lb{2.20} \\
\ti\calC_{m\ell} &= \langle x_m, zx_\ell\rangle \lb{2.21}
\end{align}
called the {\it CMV matrix\/} and the {\it alternate CMV matrix}, respectively.
By \eqref{2.17}, the unitarity of $\calC$ and $\bar z = z^{-1}$, we see
\begin{equation} \lb{2.22}
\ti\calC_{mk} = \calC_{km}
\end{equation}
that is, $\calC$ and $\ti\calC$ are transposes of each other. We will see shortly that
$\calC$ is five-diagonal, but this follows now by noting that both $z$ and $z^{-1}$ map
$[\chi_0, \dots, \chi_k]$ into $[\chi_0, \dots, \chi_{k+2}]$.

CMV \cite{CMV}, Ammar--Gragg--Reichel \cite{Athens}, Bunse-Gerstner and Elsner \cite{BGE},
and Watkins \cite{Watk} also discussed the important factorization $\calC=\calL\calM$ as follows:
\begin{equation} \lb{2.23x}
\calL_{mk} = \langle \chi_m, zx_k\rangle \qquad
\calM_{mk} = \langle x_m, \chi_k\rangle
\end{equation}
Since $\{x_k\}_{k=1}^\infty$ is a basis, $\langle f,g\rangle = \sum_{k=0}^\infty
\langle f,x_k\rangle \langle x_k,g\rangle$, and thus
\begin{equation} \lb{2.23}
\calC=\calL\calM \qquad \ti\calC=\calM\calL
\end{equation}
The point of this factorization is that $\calL$ and $\calM$ have a simpler structure than
$\calC$. Indeed, $\calL$ is a direct sum of $2\times 2$ blocks and $\calM$ of a single
$1\times 1$ block and then $2\times 2$ blocks.

One can (and we will) see this based on calculations, but it is worth seeing why it is true
in terms of the structure of the CMV and alternate CMV basis. Notice that $\chi_{2n-1}$
and $\chi_{2n}$ span the two-dimensional space $[1,z,z^{-1}, \dots, z^n, z^{-n}]\cap
[1,z,z^{-1}, \dots, z^{n-1},z^{-n+1}]^\perp$ and so do $x_{2n-1}$ and $x_{2n}$.
This shows that $\calM$ is a direct sum of $\bdone_{1\times 1}$ and $2\times 2$ matrices.
Similarly, $\chi_{2n}$ and $\chi_{2n+1}$ span $[1,\dots, z^{-n}, z^{n+1}]\cap
[1,\dots, z^{-n+1}, z^n]^\perp$, as do $zx_{2n}$ and $zx_{2n+1}$ (even for
$n=0$). Thus $\calL$ has a $2\times 2$ block structure.

In fact, we can use Szeg\H{o} recursion in the form \eqref{2.5}, \eqref{2.6}, \eqref{2.10},
\eqref{2.11} to find the $2\times 2$ matrices explicitly. For example, taking \eqref{2.11}
for $n=2m-1$, we get
\[
\varphi_{2m}^* = \rho_{2n-1} \varphi_{2m-1}^* - \alpha_{2m-1} \varphi_{2m}
\]
and multiplying by $z^{-m}$ yields (by \eqref{2.19}/\eqref{2.19}),
\[
\chi_{2m} = -\alpha_{2m-1} x_{2m} + \rho_{2m-1} x_{2m-1}
\]

This plus similar calculations imply

\begin{theorem}[\cite{Athens,BGE,CMV,Watk}] \lb{T2.1}  Let
\begin{equation} \lb{2.24}
\Theta(\alpha) = \begin{pmatrix}
\bar\alpha & \rho \\
\rho & -\alpha \end{pmatrix}
\end{equation}
Then $\calC=\calL\calM$ and
\begin{align}
\calL &= \Theta(\alpha_0) \oplus \Theta(\alpha_2) \oplus \Theta(\alpha_4) \oplus\cdots
\oplus\Theta (\alpha_{2m}) \oplus\cdots \lb{2.25} \\
\intertext{and}
\calM &=\bdone_{1\times 1} \oplus\Theta (\alpha_1) \oplus\Theta (\alpha_3) \oplus\cdots
\oplus\Theta (\alpha_{2m+1})\oplus\cdots  \lb{2.26}
\end{align}
\end{theorem}

Doing the multiplication yields
\begin{equation} \lb{2.27}
\calC = \begin{pmatrix}
{}& \bar\alpha_0 & \bar\alpha_1 \rho_0 & \rho_1 \rho_0 & 0 & 0 & \dots & {} \\
{}& \rho_0 & -\bar\alpha_1 \alpha_0 & -\rho_1 \alpha_0 & 0 & 0 & \dots & {} \\
{}& 0 & \bar\alpha_2\rho_1 & -\bar\alpha_2 \alpha_1 & \bar\alpha_3 \rho_2 & \rho_3 \rho_2 & \dots & {} \\
{}& 0 & \rho_2 \rho_1 & -\rho_2 \alpha_1 & -\bar\alpha_3 \alpha_2 & -\rho_3 \alpha_2 & \dots & {} \\
{}& 0 & 0 & 0 & \bar\alpha_4 \rho_3 & -\bar\alpha_4 \alpha_3 & \dots & {} \\
{}& \dots & \dots & \dots & \dots & \dots & \dots & {}
\end{pmatrix}
\end{equation}

$\calC$ is five-diagonal, that is, only nonzero in those diagonals $\calC_{k\, k+j}$ with
$j=0,\pm 1, \pm 2$. Notice that half of the elements with $j=\pm 2$ are zero, so it is only
``barely" five-diagonal --- and it cannot be tridiagonal or even four-diagonal since

\begin{proposition}[\cite{CMV2}]\lb{P2.1A} If $\{A_{jk}\}_{1\leq j,\, k<\infty}$ is a
semi-infinite unitary matrix and
\[
k-j \notin \{-1, \dots, n\} \Rightarrow A_{jk}=0
\]
then $A$ is a direct sum of finite blocks of size at most $n+1$.
\end{proposition}

This was proven for $n=1$ in \cite{BHJ} and conjectured for $n=2$ in a draft of
\cite{OPUC1} before motivating \cite{CMV2}.

While our construction has been for $\alpha_n$'s which come from a $d\mu$ and, in
particular, which obey
\begin{equation} \lb{2.28x}
\abs{\alpha_j (d\mu)} <1
\end{equation}
$\Theta$ defines a unitary so long as $\abs{\alpha_n}\leq 1$. We thus define
a {\it CMV matrix\/} to be a matrix of the form \eqref{2.23}--\eqref{2.27} for any
$\{\alpha_n\}_{n=0}^\infty$ with $\abs{\alpha_n} \leq 1$. If $\abs{\alpha_n}<1$ for all
$n$, we call $\calC$ a {\it proper CMV matrix}, and if $\abs{\alpha_n}=1$ for some
$n$, we call it an {\it improper CMV matrix}.

To state the analog of Stone's selfadjoint cyclic model theorem, we need another
definition. A {\it cyclic unitary model\/} is a unitary operator, $U$\!, on a (separable)
Hilbert space, $\calH$, with a distinguished unit vector, $v_0$, which is cyclic,
that is, finite linear combinations of $\{U^n v_0\}_{n=-\infty}^\infty$ are dense
in $\calH$. We call the model {\it proper\/} if $\dim(\calH)=\infty$ and
{\it improper\/} if $\dim(\calH)<\infty$. It is easy to see that the model is
improper if and only if $P(U)=0$ for some polynomial, $P$, which can be taken
to have degree $\dim(\calH)-1$. Two cyclic unitary models, $(\calH, U,v_0)$
and $(\wti\calH, \wti U, \ti v_0)$, are called equivalent if and only if there
is a unitary $W$ from $\calH$ onto $\wti\calH$ so that
\begin{equation} \lb{2.28}
Wv_0 =\ti v_0 \qquad WUW^{-1} =\wti U
\end{equation}

\begin{theorem}\lb{T2.2} There is a one-one correspondence between proper cyclic
unitary models and proper CMV matrices, $\calC$, in that $\delta_0 = (1,0,0,\dots)^t$ is
cyclic for any such $\calC$ and every equivalence class contains exactly one proper
CMV model: $(\ell^2,\calC, \delta_0)$.
\end{theorem}

\begin{remarks} 1. There is behind this a four-fold equivalence:
\begin{SL}
\item[(i)] equivalence classes of proper cyclic unitary models
\item[(ii)] nontrivial probability measures on $\partial\bbD$
\item[(iii)] Verblunsky coefficients $\{\alpha_n (d\mu)\}_{n=0}^\infty$ in
$\bbD^\infty$
\item[(iv)] proper CMV matrices.
\end{SL}
The spectral theorem sets up a one-one correspondence between (i) and (ii), while
the definition of CMV matrices between (iii) and (iv). Szeg\H{o} recursion sets up
a map from $d\mu$ to $\{\alpha_n(d\mu)\}_{n=1}^\infty$. As we will show, each $(\ell^2,
\calC, \delta_0)$ is a cyclic model, so the key remaining fact is the uniqueness.

\smallskip
2. A corollary of this is Verblunsky's theorem (also called ``Favard's theorem for
the unit circle") that each $\{\alpha_n\}_{n=0}^\infty\in\bbD$ is the Verblunsky
coefficient for some $d\mu$. See \cite{OPUC1,1Foot} for further discussion and other
proofs.
\end{remarks}

\begin{proof} As explained in Remark~1, we need only prove that any proper CMV matrix
has $\delta_0$ as a cyclic vector, and that if $\{\alpha_n^{(0)}\}_{n=0}^\infty$ are
the Verblunsky coefficients for $\calC$ and $d\mu$ the spectral measure for $\delta_0$, then
\begin{equation} \lb{2.29}
\alpha_n (d\mu) =\alpha_n^{(0)}
\end{equation}

Let $\delta_n$ be the unit vector in $\ell^2$ with coefficient $1$ in place
$n$ and $0$ elsewhere; index labelling for our vectors starts at $0$.
By direct calculations using the $\calL\calM$ representation,
\begin{align}
\calC^{n+1} \delta_0 -\rho_0^{(0)} \rho_1^{(0)} \dots \rho_{2n}^{(0)}
\delta_{2n+1} &\in [\delta_0, \dots, \delta_{2n}] \lb{2.30} \\
(\calC^*)^n \delta_0 - \rho_0^{(0)} \dots \rho_{2n-1}^{(0)}\delta_{2n} &\in
[\delta_0, \dots, \delta_{2n-1}] \lb{2.31}
\end{align}
It follows that $\delta_0$ is cyclic and
\begin{equation} \lb{2.32}
\chi_n(d\mu) =W\delta_n
\end{equation}
where $W$ is the spectral representation from $\ell^2$ to $L^2 (\partial\bbD,d\mu)$.
\eqref{2.32} follows from \eqref{2.30}--\eqref{2.31}, induction, and the Gram--Schmidt
definition of $\chi$.

By \eqref{2.32} and
\begin{alignat*}{2}
\langle \delta_0,\calC\delta_0\rangle &=\bar\alpha_0^{(0)} &\qquad
\langle \chi_0, z\chi_0\rangle &= \alpha_0(d\mu) \\
\langle \delta_{2n-2}, \calC\delta_{2n-1}\rangle &= \bar\alpha_{2n-1}^{(0)} \rho_{2n-2}^{(0)}
& \qquad \langle \chi_{2n-2}, z\chi_{2n-1}\rangle &=\bar\alpha_{2n-1} \rho_{2n-2} \\
\langle \delta_{2n}, \calC\delta_{2n-1}\rangle &= \bar\alpha_{2n}^{(0)} \rho_{2n-1}^{(0)}
& \qquad \langle \chi_{2n}, z\chi_{2n-1}\rangle &= \bar\alpha_{2n}(d\mu) \rho_{2n-1}(d\mu)
\end{alignat*}
we obtain \eqref{2.29} by induction.
\end{proof}

\section{Cutoff, Finite, Two-Sided, Periodic, and Floquet CMV Matrices} \lb{s3}

In this section, we will discuss various matrices constructed from or related to
CMV matrices. Some are finite, and in that case, we will also discuss the associated
characteristic polynomial which turns out to be equal or related to the basic ordinary
or Laurent polynomials of OPUC: the monic orthogonal and paraorthogonal polynomials
and the discriminant. The basic objects we will discuss are:

\begin{SL}
\item[(i)] {\it Cutoff CMV matrices}, that is, $\ti\pi_n \calC\ti\pi_n$ where
$\ti\pi_n$ is projection onto the span of the first $n$ of $1,z,z^{-1},\dots$.

\item[(ii)] {\it Finite CMV matrices}, the upper $n\times n$ block of an improper
CMV matrix with $\alpha_{n-1}\in\partial\bbD$.

\item[(iii)] {\it Two-sided CMV matrices\/} defined for $\{\alpha_n\}_{n=-\infty}^\infty$
via extending $\calL$ and $\calM$ in the obvious way to a two-sided form.

\item[(iv)] {\it Periodic CMV matrices}. The special case of two-sided CMV matrices
when $\alpha_{n+p}=\alpha_n$ for some $p$.

\item[(v)] {\it Floquet CMV matrices}. Periodic CMV matrices have a direct
integral decomposition whose fibers are $p\times p$ matrices that are finite
CMV matrices with a few changed matrix elements.
\end{SL}

\smallskip
\noindent{\bf Cutoff CMV matrices}. \  A {\it cutoff CMV matrix\/} is the restriction
of a proper CMV matrix to the upper $n\times n$ block, that is, top $n$ rows and leftmost
$n$ columns. We use $\calC^{(n)}$ to denote the cutoff matrix associated to $\calC$.
A glance at \eqref{2.27} shows that $\calC^{(n)}$ depends on $\{\alpha_j\}_{j=0}^{n-1}$.
Here is a key fact:

\begin{proposition}\lb{P3.1} Let $\Phi_n(z)$ be the monic orthogonal polynomial
associated to $\calC$ {\rm{(}}i.e., $\Phi_n = \kappa_n^{-1}\varphi_n${\rm{)}}. Then
\begin{equation} \lb{3.1}
\Phi_n(z) = \det (z\bdone -\calC^{(n)})
\end{equation}
\end{proposition}

\begin{proof} If $\pi_n$ is the projection onto $[1,\dots, z^{n-1}]$ and $\ti\pi_n$
on the span of the first $n$ of $1,z,z^{-1}, \dots$, then $\pi_n$ and $\ti\pi_n$
are unitarily equivalent under a power of $z$. So if $M_z f=zf$, then $\pi_n M_z
\pi_n$ and
\begin{equation} \lb{3.1a}
\ti\pi_n M_z \ti\pi_n \equiv \calC^{(n)}
\end{equation}
are unitarily equivalent, and thus, \eqref{3.1} is equivalent to
\begin{equation} \lb{3.2}
\Phi_n(z) =\det (z\bdone -\pi_n M_z \pi_n)
\end{equation}

Let $z_j$ be a zero of $\Phi_n$ of multiplicity $k_j$ and let $P_j(z) =
\Phi_n(z)/(z-z_j)^{k_j}$. Then with $A=\pi_n M_z \pi_n$, we have
\[
(A-z_j)^{k_j} P_j =0 \qquad (A-z_j)^{k_j-1} P_j\neq 0
\]
Thus, as $z_j$ runs through the distinct zeros, $\{(A-z_j)^\ell P_j \mid
\ell=0,1,\dots, k_j-1\}$ gives us a Jordan basis in which $A$ has a
$k_j\times k_j$ block for each $z_j$ of the form
\[
\begin{pmatrix}
z_j & 1 & \dots & \dots & 0 \\
0 & z_j & 1 & \dots & 0 \\
\hdotsfor{5} \\
0 & 0 & \dots & \dots & 1 \\
0 & 0 & \dots & \dots & z_j
\end{pmatrix}
\]
and thus
\[
\det(z-A) = \prod (z-z_j)^{k_j} = \Phi_n (z)
\qedhere
\]
\end{proof}

\begin{corollary} \lb{C3.2} The zeros of $\Phi_n(z)$ lie in $\bbD$.
\end{corollary}

\begin{remark} See Section~1.7 of \cite{OPUC1} for six other proofs of this theorem.
\end{remark}

\begin{proof} Let $A=\pi_n M_z\pi_n$. Then $\|A\|\leq 1$, so obviously, eigenvalues lie
in $\ol{\bbD}$. If $A\eta =z_0\eta$ with $\eta\in\ran \pi_n$ and $z_0\in\partial\bbD$, then
$\|A\eta\|=\|\eta\|$, so $\pi_n z\eta = z\eta$ and thus as polynomials $(z-z_0)\eta =0$.
Since the polynomials are a division ring, $\eta =0$, that is, there are no eigenvalues
on $\partial\bbD$.
\end{proof}

To classify cutoff CMV matrices, we need to understand how $\Theta(\alpha)$ arises
from a $2\times 2$ change of basis.

\begin{lemma}\lb{L3.2A} Let $f,g$ be two independent unit vectors with
\begin{equation} \lb{3.3}
\langle f,g \rangle =\alpha
\end{equation}
Let $\varphi_1, \varphi_2$ be the result of applying Gram--Schmidt to $f,g$ and
$\psi_1,\psi_2$ to $g,f$. Let $M$ be the matrix of change of basis
\begin{align}
\varphi_1 &= m_{11} \psi_1 + m_{12} \psi_2  \lb{3.4a} \\
\varphi_2 &= m_{21} \psi_1 + m_{22} \psi_2  \lb{3.4b}
\end{align}
Then $\alpha\in\bbD$ and $M=\Theta(\alpha)$.
\end{lemma}

\begin{proof} $\abs{\alpha}<1$ by the Schwarz inequality and the independence of
$f$ and $g$. Note that
\begin{align*}
\|g-\langle f,g\rangle f\|^2 &= \|g\|^2 + \abs{\langle f,g\rangle}^2 \|f\|^2
- 2\Real [\langle f,g\rangle \langle g,f\rangle ] \\
&= 1-\abs{\alpha}^2 \equiv \rho^2
\end{align*}
so
\begin{align}
\varphi_2 &= \rho^{-1} (g-\alpha f) \lb{3.5} \\
\psi_2 &= \rho^{-1} (f-\bar\alpha g) \lb{3.6}
\end{align}

From this, a direct calculation shows that
\begin{align*}
m_{11} &=\langle \psi_1,\varphi_1\rangle = \langle g,f\rangle =\bar\alpha \\
m_{12} &=\langle \psi_2,\varphi_1\rangle = \rho^{-1} (1-\abs{\alpha}^2) =\rho \\
m_{21} &=\langle \psi_1,\varphi_2\rangle = \rho^{-1} (1-\abs{\alpha}^2) =\rho \\
m_{22} &=\langle \psi_2,\varphi_2\rangle = \rho^{-2} (\alpha + \alpha \abs{\alpha}^2
-2\alpha)= -\alpha
\qedhere
\end{align*}
\end{proof}

\begin{remark} One can use this lemma to deduce the form of $\calL$ and $\calM$
in the $\calL\calM$ factorization without recourse to Szeg\H{o} recursion, and
then use their form to deduce the Szeg\H{o} recursion. This is precisely
Watkins' approach \cite{Watk} to the factorization and Szeg\H{o} recursion.
\end{remark}

Given any matrix $A$ and vector $\delta_0$, define
\begin{equation} \lb{3.7}
V_k = \Span (\delta_0, A\delta_0, A^* \delta_0, \dots, A^k \delta_0,
(A^*)^k\delta_0)
\end{equation}
which has dimension $2k+1$ if and only if the vectors are independent.

Here is a complete classification of cutoff CMV matrices analogous to
Theorem~\ref{T2.2}:

\begin{theorem}\lb{T3.3} Let $A$ be an $n\times n$ cutoff CMV matrix with
$\delta_0 = (1,0,\dots, 0)^t$. Then:
\begin{SL}
\item[{\rm{(i)}}] If $n=2k+1$, $V_k$ has dimension $n$. If $n=2k$, then
$\Span [V_{k-1} \cup \{A^k \delta_0\}]$ has dimension $n$.
\item[{\rm{(ii)}}] If $n=2k+1$, $A^*$ is an isometry on $\Span (V_{k-1}
\cup \{(A)^k\delta_0\})$ and $A$ is an isometry on $\Span (V_{k-1}\cup
\{(A^*)^k \delta_0\})$. If $n=2k$, $A^*$ is an isometry on $\Span
(V_{k-2}\cup \{A^{k-1} \delta_0, A^k\delta_0\})$ and $A$ is an isometry
on $V_{k-1}$.
\item[{\rm{(iii)}}] $\|A\|\leq 1$.
\item[{\rm{(iv)}}] $A$ is not unitary.
\end{SL}

Conversely, if $A_0,\delta_0$ are a pair of an $n\times n$ matrix and
vector $\delta_0$ obeying {\rm{(i)--(iv)}}, then there is a basis in which
$\delta_0 = (1,0, \dots, 0)^t$ and $A$ is a cutoff CMV matrix.

$(A,\delta_0)$ determine the Verblunsky coefficients $(\alpha_0, \dots,
\alpha_{n-1})$ uniquely. In particular, two cutoff CMV matrices with distinct
$\{\alpha_j\}_{j=0}^n$ are not unitarily equivalent by a unitary preserving
$\delta_0$.
\end{theorem}

\begin{proof} Suppose first that $A$ is a cutoff CMV matrix, that is,
$A=\ti\pi_n \calC\ti\pi_n$. By definition of $\ti\pi_n$,
\[
\ti\pi_n \calC^j \delta_0 = \calC^j\delta_0 \qquad
\ti\pi_n (\calC^*)^\ell \delta_0 = (\calC^*)^\ell \delta_0
\]
for $j=0,1,\dots,k$ and $\ell =0,1,\dots, k$ (resp., $k-1$) if $n=2k+1$
(resp., $2k)$. It follows that for those values of $j$ and $\ell$,
\begin{equation} \lb{3.8}
\calC^j\delta_0 = A^j\delta_0 \qquad
(\calC^*)^\ell \delta_0 = (A^*)^\ell \delta_0
\end{equation}
so that (i) holds.

This also shows that $A^* A^j\delta_0 = A^{j-1} \delta_0$ for $j=1, \dots, k$,
and from this and \eqref{3.8}, it follows that $A^*$ is unitary on
$\Span \{A^j \delta_0\}_{j=0}^k\cup \{(A^*)^\ell
\delta_0\}_{\ell=0}^{k-1 \text{ (or $k-2$)}}$. Similarly, we get
the unitarity result for $A$.

(iii) is obvious since $\|\ti\pi_n\| = \|\calC\|=1$ and (iv) follows since there
is a vector $\varphi$ in $\ran (\ti\pi_n)$ with $\ti\pi_n \calC\varphi\neq \calC
\varphi$. This completes the proof of the first paragraph of the theorem.

As a preliminary to the converse, we note that $\ti\pi_n$ commutes with either
$\calL$ or $\calM$, so a finite CMV matrix has the form $\calL_n\calM_n$ where
if $n=2k+1$ is odd ($\bdone=1\times 1$ identity matrix),
\begin{align}
\calL_n &= \Theta (\alpha_0)\oplus\cdots\oplus\Theta (\alpha_{2k-2})
\oplus \alpha_{2k} \bdone \lb{3.9a} \\
\calM_n &= \bdone \oplus \Theta (\alpha_1)\oplus\cdots\oplus \Theta (\alpha_{2k-1})
\lb{3.9b}
\end{align}
and if $n=2k$ is even,
\begin{align}
\calL_n &= \Theta (\alpha_0)\oplus\cdots\oplus \Theta (\alpha_{2k-2}) \lb{3.10a} \\
\calM_n &= \bdone \oplus\Theta (\alpha_1)\oplus\cdots\oplus\Theta (\alpha_{2k-3})
\oplus \alpha_{2k-1}\bdone \lb{3.10b}
\end{align}
We will prove that when $A$ obeys (i)--(iv), then $A$ has an $\calL_n\calM_n$
factorization with parameter $\alpha_j$ given intrinsically by $A$. This will
not only prove the converse but the uniqueness of the map from
$\{\alpha_j\}_{j=0}^{N-1}$ to cutoff CMV matrices, and so it will complete
the proof of the theorem.

We first consider the case $n=2k+1$ odd. Define $\chi_\ell$ to be the basis obtained
by Gram--Schmidt on $\delta_0, A\delta_0, A^*\delta_0,\dots, A^k\delta_0,
(A^*)^k \delta_0$ (this is possible because (i) implies these vectors are
linearly independent) and define $x_\ell$ to be the result of Gram--Schmidt
on $\delta_0, A^*\delta_0, A\delta_0, \dots, (A^*)^k \delta_0, A^k\delta_0$.
Then if $A$ is written in $\chi_\ell$ basis,
\begin{equation} \lb{3.11a}
A= \calL \calM
\end{equation}
where
\begin{equation} \lb{3.11}
\calM_{k\ell} =\langle x_k, \chi_\ell\rangle \qquad
\calL_{k\ell} = \langle \chi_k, Ax_\ell\rangle
\end{equation}
We need to show that $\calL,\calM$ have the form \eqref{3.9a}/\eqref{3.9b}.

If $P_m$ is the projection to the orthogonal complement of $V_{m-1}$ and
$f=P_{m-1}(A^*)^m\delta_0/\|P_{m-1}(A^*)^m \delta_0\|$ and $g=P_{m-1} A^m\delta_0/
\|P_{m-1} A^m \delta_0\|$, then $\{\chi_\ell, x_\ell\}_{\ell=m, m+1}$ are
given by Lemma~\ref{L3.2A}. So $M$ has the form $1\oplus \Theta (\alpha_1)
\oplus\cdots\oplus\Theta (\alpha_{2k-1})$ as required.

Let $W_\ell$ be the projection onto the span of the first $2\ell$ of $\delta_0,
A\delta_0, A^*\delta_0, A^2\delta_0, \dots$ and $\wti W_\ell$ the span of the
first $2\ell$ of $\delta_0, A^*\delta_0, A\delta_0, (A^*)^2\delta_0, \dots$.
By hypothesis (ii), $A$ is an isometry on $\wti W_1, \wti W_2, \dots,
\wti W_k$, and by the same hypothesis, $AA^*\varphi=\varphi$ for $\varphi
=\delta_0, A^*\delta_0, \dots, (A^*)^k\delta_0$. So it follows that $A$ maps
$\wti W_\ell$ to $W_\ell$ for $\ell=1, \dots, k$. Thus, by Lemma~\ref{L3.2A},
the $2k\times 2k$ upper block of $L$ is $\Theta(\alpha_0) \oplus\Theta
(\alpha_2)\oplus\cdots\oplus\Theta (\alpha_{2k-2})$. Since $A$ and $A^*$ are
contractions, $L$ must have $0$'s in the bottom and rightmost column, except
for the lower corner. That corner value, call it $\alpha_{2k}$, must have
$\abs{\alpha_{2k}}\leq 1$ by (iii) and $\abs{\alpha_{2k}}<1$ by (iv).
Thus, we have the required $\calL\calM$ factorization if $n=2k+1$.

Now let $\ell=2k$ be even. Define $\chi_\ell$ as before, but define
$\ti x_\ell$ by Gram--Schmidt on $A\delta_0, \delta_0, A^2\delta_0,
A^*\delta_0, \dots, A^k\delta_0, (A^*)^{k-1}\delta_0$. Then $A$ written
in $\chi_\ell$ basis has the form \eqref{3.11a} where
\begin{equation} \lb{3.12}
\calM_{k\ell}=\langle\ti x_k, A\chi_\ell\rangle \qquad
\calL_{k\ell} =\langle \chi_k, \ti x_\ell\rangle
\end{equation}
We need to show that $\calL,\calM$ have the form \eqref{3.10a}/\eqref{3.10b}.

Define $W_\ell$ to be the span of $\delta_0, A\delta_0, A^*\delta_0,
\dots, A^\ell \delta_0, (A^\ell)^*\delta_0$ and $\wti W_\ell$ the span
of $A\delta_0, \delta_0, \dots, (A^*)^{\ell-1}\delta_0, A^{\ell+1}\delta_0$.
As above, $A$ is an isometry of $W_\ell$ to $\wti W_\ell$, so $M$ has the form
$\bdone\oplus\Theta(\alpha_1)\oplus\cdots\oplus\Theta (\alpha_{2k-3})\oplus
\alpha_{2k-1}\bdone$ where $\abs{\alpha_{2k-1}}<1$ by condition (iv).
Similarly, $L$ has a $\Theta(\alpha_0)\oplus\cdots\oplus\Theta(\alpha_{2k-2})$
block structure. This proves \eqref{3.10a}/\eqref{3.10b} and completes the
case $n=2k$.
\end{proof}

\begin{remark} This theorem sets up a one-one correspondence between
$\{\alpha_j\}_{j=0}^{n-1}\in\bbD^n$ and cutoff CMV matrices.
\end{remark}
\medskip

\noindent{\bf Finite CMV matrices.} \ As discussed in Section~\ref{s2}, $\calC$, originally
defined for $\abs{\alpha_j}<1$, has an extension to $\abs{\alpha_j}\leq 1$ via the
$\Theta$ formula for $\calL,\calM$. Since $\abs{\alpha_{j_0}}=1$ implies $\rho_{j_0}=0$
and $\Theta(\alpha_{j_0})$ is diagonal, if $\abs{\alpha_{j_0}}=1$, $\calC(\{\alpha_j\})$
leaves $\bbC^{j_0+1}$ (i.e., vectors $\varphi$ with $\varphi_k=0$ if $k\geq j_0 +1$)
invariant and $\calC\restriction\bbC^{j_0+1}$ is a $(j_0+1) \times (j_0+1)$
unitary matrix. If $\abs{\alpha_0}, \dots,\abs{\alpha_{n-2}}<1$ and
$\abs{\alpha_{n-1}}=1$, the corresponding $n\times n$ matrix is called a
{\it finite $n\times n$ CMV matrix}, $\calC_n (\{\alpha_0, \dots,
\alpha_{n-2}, \alpha_{n-1}\})$. $\calC_n$ has the form $\calL_n \calM_n$ where
\eqref{3.9a}/\eqref{3.9b} or \eqref{3.10a}/\eqref{3.10b} hold, and now
$\alpha_{n-1}\in\partial\bbD$.

If $U$ is an $n\times n$ matrix and $\delta_0$ a cyclic vector in the sense that
$\{U^m \delta_0\}_{m=-\infty}^\infty$ is total, $\delta_0$ cannot be orthogonal
to any eigenvector. So $U$ has to have $n$ distinct eigenvalues $\{\lambda_j\}_{j=1}^n$
and the eigenvectors $\{\psi_j\}_{j=1}^n$ obey $\abs{\langle \delta_0, \psi_j\rangle}^2
= a_j\neq 0$. The unitary invariants of the pair $(U,\delta_0)$ are the spectral
measures $\sum_{j=1}^n a_j \delta_{\lambda_j}$ where $\{\lambda_j\}_{j=1}^n$ are
arbitrary distinct points and the $a_j >0$ have the single restriction $\sum_{j=1}^n
a_j=1$. Thus, the number of real parameters is $n+(n-1)=2n-1$. The number of free
parameters in an $n\times n$ finite CMV matrix is $n-1$ complex numbers in $\bbD$
and one in $\partial\bbD$, that is, $2(n-1)+1=2n-1$. This suggests that

\begin{theorem}\lb{T3.4} There is a one-one correspondence between unitary
equivalence classes of $n\times n$ unitary matrices with a cyclic vector and
finite CMV matrices in that each equivalence class contains one CMV matrix
{\rm{(}}fixed by $\delta_0 = (1,0,\dots, 0)^t${\rm{)}} and two CMV matrices
with distinct parameters are not unitarily equivalent by a unitary fixing
$(1,0,\dots, 0)^t$.
\end{theorem}

The proof is identical to the proof of Theorem~\ref{T3.3} except that $A$ nonunitary
is replaced by $A$ unitary so $\abs{\alpha_{n-1}}=1$. As noted in the discussion
after Lemma~\ref{L3.2A}, this approach is close to that in Watkins \cite{Watk}.
This theorem is related to results in Ammar--Gragg--Reichel \cite{Athens} and
Killip--Nenciu \cite{KNnew}. The latter talk about matrices with CMV shape having
the CMV form.

Instead of the cutoff CMV matrix, $\pi_n M_z \pi_n$, one can look at
$\widehat\pi_n M_z \widehat\pi_n$ where $\widehat\pi_n$ is a not necessary selfadjoint
projection. CMV \cite{CMVppt} have shown that finite CMV matrices have this
form and that they are the only normal operators of this form.

\medskip
\noindent{\bf Two-sided CMV matrices.} \ In a sense, CMV matrices are two-sided. For example,
if $\alpha_n\equiv 0$, $\calC$ is unitarily equivalent to a two-sided shift since
$\calC^k\delta_0 = \delta_{2k-1}$ and $\calC^{-k}\delta_0 = \delta_{2k}$. However,
as structures, the matrix is semi-infinite and there is a cyclic vector which is
often not true for two-sided matrices. Thus, there is an extension to ``two-sided"
examples.

Let $\{\alpha_n\}_{n=-\infty}^\infty$ be a two-sided sequence of numbers in $\ol{\bbD}$.
Let $\calH=\ell^2(\bbZ)$, that is, two-sided sequences $\{u_n\}_{n=-\infty}^\infty$
with $\sum_{n=-\infty}^\infty \abs{u_n}^2 <\infty$. Let $\Theta_j(\beta)$ be
$\Theta(\beta)$ acting on the two indices $j$ and $j+1$. Define
\begin{equation} \lb{3.13}
\calE (\{\alpha_j\}_{j=-\infty}^\infty) = \ti\calL (\{\alpha_j\}_{j=-\infty}^\infty)
\wti\calM (\{\alpha_j\}_{j=-\infty}^\infty)
\end{equation}
where
\begin{align*}
\wti\calM &= \bigoplus_{j=-\infty}^\infty \Theta_{2j-1} (\alpha_{2j-1}) \\
\wti\calL &= \bigoplus_{j=-\infty}^\infty \Theta_{2j} (\alpha_{2j})
\end{align*}
$\calE$ is called the {\it extended CMV matrix}.

The extended CMV matrix was introduced in \cite{OPUC1}. Earlier,
Bourget, Howland, and Joye \cite{BHJ} had considered some doubly infinite five-diagonal
matrices which factor into a product of two direct sums of $2\times 2$ matrices,
but the $2\times 2$ blocks were general unitaries rather than $\Theta$'s.

While $\calE$ is natural and important for the periodic case, we will also see
that it arises in the theory of essential spectrum of $\calC$ (see Section~\ref{s8}).

One reason for the name ``extended CMV matrix" is:

\begin{proposition}\lb{P3.5} If $\alpha_{-1}=-1$, then $\calE$ is a direct sum
on $\ell^2 (-\infty, -1] \oplus\ell^2 [0,\infty)$ and $\calE\restriction \ell^2
[0,\infty)$ is the CMV matrix $\calC(\{\alpha_j\}_{j=0}^\infty)$. Moreover,
$\calE\restriction \ell^2 (-\infty, -1]$ is unitarily equivalent to $\calC
(\{\bar\alpha_{-j-2}\}_{j=0}^\infty)$.
\end{proposition}

\begin{remark} $\ell^2 [0,\infty)$ means those $u\in\ell^2(\bbZ)$ with $u_n=0$
if $n<0$ and $\ell^2 (-\infty, -1]$ those with $u_n =0$ if $n>-1$.
\end{remark}

\begin{proof} $\Theta(-1) =\left(\begin{smallmatrix} -1 & 0 \\ 0 & 1\end{smallmatrix}
\right)$, so both $\ti\calL$ and $\wti\calM$ leave $\ell^2 [0,\infty)$ and $\ell^2
(-\infty,-1]$ invariant. Thus, $\calE$ does.

$\wti\calM\restriction\ell^2 [0,\infty) =\calM$ and $\ti\calL\restriction\ell^2
[0,\infty) =\calL$, so $\calE\restriction\ell^2 [0,\infty)$ is $\calC
(\{\alpha_j\}_{j=0}^\infty)$.

For the restriction to $\ell^2 (-\infty, -1]$, note first that $\left(
\begin{smallmatrix} 0 & 1 \\ 1 & 0\end{smallmatrix}\right) \Theta(\alpha)
\left( \begin{smallmatrix} 0 & 1 \\ 1 & 0\end{smallmatrix}\right) =
\Theta (-\bar\alpha)$. Thus, by labeling the basis backwards, $\calE$ is
unitarily equivalent to something that looks very much like $\calC
(\{-\bar\alpha_{-j-2}\}_{j=0}^\infty)$ except $\calM$ starts with $-1$,
not $1$. By the discussion in Section~\ref{s5}, there is a unitary that flips the
spin of this $-1$ and all the $\alpha_j$'s.
\end{proof}

Changing $\alpha_{-1}$ from its value to $\alpha_{-1}=-1$ is a perturbation
of rank at most two, so by the Kato--Rosenblum theorem \cite{RS3}, the
a.c.\ spectrum of $\calE$ is that of a direct sum of two $\calC$'s. Since
these a.c.\ spectra are only restricted by simplicity, we see that the
a.c.\ spectrum of $\calE$ has multiplicity at most $2$, but is otherwise
arbitrary: it can be partially multiplicity $0$, partially $1$, and partially
$2$. In particular, $\calE$ may not have a cyclic vector.

It is a theorem of Simon \cite{S293} that the singular spectrum of $\calE$ is simple.
This is an analog of a theorem of Kac \cite{Kac62,Kac} and Gilbert \cite{Gil,Gil2} for
Schr\"odinger operators.

\medskip
\noindent{\bf Periodic CMV matrices.} \ If $\{\alpha_j\}_{j=0}^\infty$ is a
sequence of Verblunsky coefficients with
\begin{equation} \lb{3.14}
\alpha_{j+p}=\alpha_j
\end{equation}
for $j\geq 0$, $p$ fixed, and $\geq 1$, then $\alpha_j$ has a unique extension
to $j\in\bbZ$ obeying \eqref{3.14}. The corresponding $\calE$ is called a
{\it periodic CMV matrix}. The theory is simpler if $p$ is even, which we
henceforth assume. As explained in \cite{OPUC2}, there are several ways to
analyze odd $p$ once one has understood even $p$.

Associated to $\{\alpha_j\}_{j=0}^{p-1}$ is a natural Laurent polynomial,
called the discriminant, $\Delta (z;\{\alpha_j\}_{j=0}^{p-1})$,
\begin{equation} \lb{3.15}
\Delta(z) =z^{-p/2} \tr (A(\alpha_{p-1}, z) A(\alpha_{p-2},z) \dots
A(\alpha_0, z))
\end{equation}
where
\begin{equation} \lb{3.16x}
A(\alpha,z) = \rho^{-1} \begin{pmatrix}
z & - \bar\alpha \\
-z\alpha & 1
\end{pmatrix}
\end{equation}

This is analyzed in Section~11.1 of \cite{OPUC2}. $\Delta(z)$ is real on $\partial
\bbD$ and has positive leading coefficient. This means $\Delta(z)$ has $p$ real
parameters. This suggests the map from $\{\alpha_j\}_{j=0}^p$ (of real dimension
$2p$) to $\Delta$ is many to $1$, with inverse images generically of dimension $p$
($=2p-p$). This is in fact true:  the inverse images are tori of dimension
$d\leq p$ (we will see what $d$ is in a moment). They are called {\it isospectral tori}.

For fixed $\{\alpha_j\}_{j=0}^{p-1}$, $\Delta^{-1} ([-2,2])$ (which is the spectrum
of $\calE$) lies in $\partial\bbD$ and is naturally $p$ closed intervals whose
endpoints are $\Delta^{-1} (\{-2,2\})$. Generically (in $\alpha)$, the intervals
are disjoint, that is, their complement (called the gaps) is $p$ nonempty open
intervals. In general, the number of open intervals in the gaps is $d$, the
dimension of the isospectral torus.

\medskip
\noindent {\bf Floquet CMV matrices.} \ If $T \colon \ell^2 (\bbZ)\to\ell^2(\bbZ)$ by
$(Tu)_n =u_{n+p}$ and $p$ is even, and if $\alpha_n =\alpha_{n+p}$, then $T\wti\calM
T^{-1}=\wti\calM$ and $T\ti\calL T^{-1}=\ti\calL$, so
\begin{equation} \lb{3.16}
T\calE T^{-1} =\calE
\end{equation}
(We will not consider odd $p$ in detail, but we note in that case, $T\wti\calM
T^{-1}=\ti\calL$ and $T\ti\calL T^{-1}=\wti\calM$ so, since $\wti\calM^t =\wti\calM$
and $\ti\calL^t =\ti\calL$ (on account of $\Theta^t=\Theta$ where ${}^t$ is transpose),
we have that $T\calE T^{-1}=\calE^t$.)

Since $T$ and $\calE$ commute, they can be ``simultaneously" diagonalized, in this
case represented on a direct integral representation. One way of making this
explicit is to define, for each $\beta\in\partial\bbD$, the space $\ell_\beta^\infty$,
the sequences $\{u_n\}_{n=-\infty}^\infty$ obeying $u_{n+p}=\beta u_n$. This is
clearly a space of dimension $p$ since $\{u_n\}_{n=-\infty}^\infty$ mapping to
$\{u_n\}_{n=0}^{p-1}$ (i.e., restriction) maps $\ell_\beta^\infty$ to $\bbC^p$.

By \eqref{3.16}, $\calE$, which maps bounded sequences to bounded sequences, maps
$\ell_\beta^\infty$ to $\ell_\beta^\infty$, and so defines a finite-dimensional
operator $\calE_p(\beta)$ under the explicit relation of $\ell_\beta^\infty$
mentioned above. One sees
\begin{equation} \lb{3.17}
\calE_p (\beta) = \calL_p \calM_p (\beta)
\end{equation}
where
\begin{align}
\calL_p &= \Theta_0 (\alpha_0) \oplus\cdots\oplus  \Theta_{p-2} (\alpha_{p-2}) \lb{3.18} \\
\calM_p(\theta) &= \Theta_1 (\alpha) \oplus\cdots\oplus \Theta_{p-3} (\alpha_{p-3})
\oplus\Theta_{p-1}^{(\beta)} (\alpha_{p-1}) \lb{3.19}
\end{align}
where $\Theta_{p-1}^{(\beta)}(\alpha)$ acts on $\delta_{p-1}$ and $\delta_0$, and in
that (ordered) basis has the form
\begin{equation} \lb{3.20}
\begin{pmatrix}
\bar\alpha & \rho\beta \\
\rho\bar\beta & -\alpha
\end{pmatrix}
\end{equation}

$\calE_p(\beta)$ is called the {\it Floquet CMV matrix}. To make precise the
connection to $\calE$, we define the unitary Fourier transform $\calF\colon \ell^2
(\bbZ) \to L^2 (\partial\bbD, \f{d\theta}{2\pi};\bbC^p)$, the set of $L^2$
functions on $\partial\bbD$ with values in $\bbC^p$ by
\begin{equation} \lb{3.21}
(\calF u)_k (\beta) = \sum_{n=-\infty}^\infty \beta^{-n} u_{k+np}
\end{equation}
Then
\begin{equation} \lb{3.22}
(\calF\calE\calF^{-1} g)(\beta) = \calE_p(\beta) g(\beta)
\end{equation}
(For details, see Section~11.2 of \cite{OPUC2}.)

Finally, we note a general relation of the eigenvalues of $\calE_p(\beta)$ and the
discriminant, $\Delta(z)$, of \eqref{3.14}. For $z_0\in\partial\bbD$ is an eigenvalue
of $\calE_p (\beta)$ if and only if there is $(u_1, u_0)^t$ so that after a $p$-step
transfer, we get $\beta (u_1, u_0)^t$, that is, if and only if $z_0^{p/2}\beta$ is an
eigenvalue of $T_p(z_0)$. This is true if and only if $z_0^{-p/2} T_p(z_0)$ has
eigenvalues $\beta$ and $\beta^{-1}$ if and only if $\Delta(z_0) =\beta + \beta^{-1}$.
It follows that
\begin{equation} \lb{3.23}
\det (z-\calE_p (\beta))= \biggl(\, \prod_{j=0}^{p-1} \rho_j\biggr)
\bigl[ z^{p/2} [\Delta(z)-\beta-\beta^{-1}]\bigr]
\end{equation}
for both sides are monic polynomials of degree $p$ and they have the same zeros.

\section{CMV Matrices for Matrix-Valued Measures} \lb{s4}

Because of applications to perturbations of periodic Jacobi and CMV matrices
\cite{DKSprep}, interest in matrix-valued measures (say, $k\times k$ matrices) has
increased. Here we will provide the CMV basis and CMV matrices in this
matrix-valued situation; these results are new here. Since adjoints of
finite-dimensional matrices enter but we want to use ${}^*$ for Szeg\H{o}
reversed polynomials, in this section we use ${}^\dag$ for matrix adjoint.

Measures which are nontrivial in a suitable sense are described by a
sequence $\{\alpha_j\}_{j=0}^\infty$ of Verblunsky coefficients that are
$k\times k$ matrices with $\|\alpha_j\|<1$.

To jump to the punch line, we will see that $\calC$ still has an $\calL\calM$
factorization, where $\Theta(\alpha)$ is the $2k\times 2k$ matrix
\begin{equation} \lb{4.1}
\Theta (\alpha) = \begin{pmatrix}
\alpha^\dag & \rho^\capL \\
\rho^\capR & -\alpha
\end{pmatrix}
\end{equation}
where
\begin{equation} \lb{4.2}
\rho^\capL = (1-\alpha^\dag\alpha)^{1/2} \qquad
\rho^\capR = (1-\alpha\alpha^\dag)^{1/2}
\end{equation}
It is an interesting calculation to check that $\Theta$ is unitary, that is,
\begin{equation} \lb{4.3}
\begin{pmatrix}
\alpha & \rho^\capR \\
\rho^\capL &-\alpha^\dag
\end{pmatrix}
\begin{pmatrix}
\alpha^\dag & \rho^\capL \\
\rho^\capR & -\alpha
\end{pmatrix}
=
\begin{pmatrix}
\bdone & 0 \\
0 & \bdone
\end{pmatrix}
\end{equation}
That $\alpha\alpha^\dag + (\rho^\capR)^2 = 1 = (\rho^\capL)^2 + \alpha^\dag\alpha$ follows
from \eqref{4.2}. That $\alpha\rho^\capL - \rho^\capR\alpha = \rho^\capL\alpha^\dag -
\alpha^\dag \rho^\capR =0$ follows by expanding the square roots in \eqref{4.2} in a Taylor
series and using $\alpha (\alpha^\dag\alpha)^m = (\alpha\alpha^\dag)^m \alpha$.

To describe the model specifically, we have a $k\times k$ matrix-valued (normalized,
positive) measure which can be described as follows: $d\mu_t(\theta)$ is a positive
scalar measure on $\partial\bbD$ and for a.e.\ $(d\mu_t (\theta))$ a matrix $A(\theta)$
obeying
\begin{equation} \lb{4.4}
A(\theta) \geq 0 \qquad \tr (A(\theta))=1
\end{equation}
We write $d\mu(\theta)=A(\theta)\, d\mu_t (\theta)$. We assume $d\mu$ is normalized
in the sense that $\int A(\theta) \, d\mu_t =\bdone$. We will consider $\calH_\capR$
to be the $k\times k$ matrix-valued functions, $f$, on $\partial\bbD$ with
\begin{equation} \lb{4.5}
\int \tr (f(\theta)^\dag A(\theta) f(\theta))\, d\mu_t(\theta) <\infty
\end{equation}
The measure $d\mu$ is called nontrivial if
\begin{equation} \lb{4.5a}
\dim [\Span \{B_\ell z^\ell\}_{\ell=0}^{n-1}] =nk^2
\end{equation}
for each $n$. Equivalently, for each $n$ and $\{B_\ell\}_{\ell=0}^{n-1}$ in
$\calL(\bbC^k)$, we have $\sum_{\ell=0}^{n-1} B_\ell z^\ell =0$ in $\calH_\capR$
implies $B_0=B_1 =\cdots = B_{n-1}=0$. Also equivalent is that $\langle\varphi,
A(\theta)\varphi\rangle \, d\mu_t(\theta)$ is nontrivial for each $\varphi\in
\bbC^k\backslash\{0\}$.

Similarly, we define $\calH_\capL$ to be $f$'s with
\begin{equation} \lb{4.5b}
\int \tr (f(\theta) A(\theta) f^\dag (\theta))\, d\mu_t (\theta) <\infty
\end{equation}
It is easy to see that nontriviality implies \eqref{4.5a} holds also in $\calH_\capL$.

We define two ``inner products," sesquilinear forms from $\calH_\capR$ and $\calH_\capL$
to $\calL (\bbC^k)$, the $k\times k$ matrices:
\begin{align}
\llangle f,g \rrangle_\capR &= \int f^\dag (\theta)\, d\mu(\theta) g(\theta)
\lb{4.6} \\
\langle\!\langle f,g \rangle\!\rangle_\capL &= \int g(\theta) \, d\mu(\theta) f^\dag (\theta)
\lb{4.7}
\end{align}
The right side of \eqref{4.6} is shorthand for
\[
\int f^\dag (\theta) A(\theta) g(\theta)\, d\mu_t (\theta)
\]
so the LHS of \eqref{4.5} is $\tr(\llangle f,g \rrangle_\capR )$. The symbols
$\capL,\capR$ (for left and right) come from
\begin{align}
\llangle f,gB \rrangle_\capR &= \llangle f,g \rrangle B \lb{4.8} \\
\llangle f,Bg \rrangle_\capL &= B \llangle f,g \rrangle \lb{4.9}
\end{align}
for constant $k\times k$ matrices, $B$.

The normalized matrix OPUC, $\varphi_n^\capR, \varphi_n^\capL$, are polynomials in
$z$ of degree $n$ with matrix coefficients with
\begin{equation} \lb{4.10}
\llangle \varphi_n^\capR, \varphi_m^\capR \rrangle_\capR = \delta_{nm}\bdone
\qquad
\llangle \varphi_n^\capL, \varphi_m^\capL \rrangle_\capR = \delta_{nm}\bdone
\end{equation}
This determines $\varphi$ uniquely up to a unitary right (resp., left) prefactor. We
will pick this prefactor by demanding
\begin{gather}
\varphi_n^{\capR,\capL}(z) = \kappa_n^{\capR,\capL} z^n + \text{lower order} \lb{4.11} \\
\kappa_{n+1}^\capL (\kappa_n^\capL)^{-1} > 0 \qquad
(\kappa_n^\capR)^{-1} \kappa_{n+1}^\capR > 0 \lb{4.12}
\end{gather}
With this choice of normalization, one has a sequence of $k\times k$ matrices,
$\{\alpha_n\}_{n=0}^\infty$, and the recursion relations
\begin{align}
z\varphi_n^\capL &= \rho_n^\capL \varphi_{n+1}^\capL
+ \alpha_n^\dag (\varphi_n^\capR)^* \lb{4.13} \\
z\varphi_n^\capR &= \varphi_{n+1}^\capR \rho_n^\capR
+ (\varphi_n^\capL)^* \alpha_n^\dag \lb{4.14} \\
(\varphi_n^\capL)^* &= (\varphi_{n+1}^\capL)^* \rho_n^\capL
+ z\varphi_n^\capR \alpha_n \lb{4.15} \\
(\varphi_n^\capR)^* &= \rho_n^\capR \varphi_{n+1}^*
+ \alpha_n z\varphi_n^\capL \lb{4.16}
\end{align}
where $\rho_n^\capR, \rho_n^\capL$ are given by \eqref{4.2} and $P_n^*(z) =z^n P_n
(1/\bar z)^\dag$. For construction of $\varphi_n^{\capL,\capR}$ and proof of
\eqref{4.13}--\eqref{4.16}, see \cite{AN84} or \cite[Section~2.13]{OPUC1} following
Delsarte et al.\ \cite{DGK78} and Geronimo \cite{Ger81}.

It will help to also have the following, which can be derived from
\eqref{4.13}--\eqref{4.16}:
\begin{align}
\varphi_{n+1}^\capL &= \rho_n^\capL z\varphi_n^\capL
- \alpha_n^\dag (\varphi_{n+1}^\capR)^* \lb{4.17} \\
(\varphi_{n+1}^\capR)^* &= \rho_n^\capR (\varphi_n^\capR)^*
- \alpha_n \varphi_{n+1}^\capL  \lb{4.18}
\end{align}

Following \eqref{2.18} and \eqref{2.19}, we define the CMV and alternate CMV
basis by
\begin{alignat}{2}
\chi_{2k}(z) &= z^{-k} (\varphi_{2k}^\capR(z))^*
\qquad && \chi_{2k-1}(z) = z^{-k+1} \varphi_{2k-1}^\capL (z) \lb{4.19} \\
x_{2k}(z) &= z^{-k} \varphi_{2k}^\capL(z)
\qquad && x_{2k-1}(z) = z^{-k} (\varphi_{2k-1}^\capR (z))^* \lb{4.20}
\end{alignat}

\begin{proposition} \lb{P4.1} $\{\chi_\ell (z)\}_{\ell=0}^\infty$ and
$\{x_\ell(z)\}_{\ell=0}^\infty$ are $\llangle \cdot \, ,\cdot\rrangle_\capL$
orthonormal, that is,
\begin{equation} \lb{4.21}
 \llangle \chi_\ell,\chi_m \rrangle_\capL=\delta_{\ell m} \qquad
\llangle x_\ell, x_m \rrangle_\capL=\delta_{\ell m}
\end{equation}
Moreover, $\chi_\ell$ is in the module span of the first $\ell$ of $1,z,
z^{-1}, \dots$ and $x_\ell$ of $1,z^{-1}, z, \dots$.
\end{proposition}

\begin{remark} By module span of $\{f_j(z)\}_{j=1}^m$ of scalar functions,
$f$, we mean elements in $\calH_\capL$ of the form $\sum_{j=1}^m B_j f_j(z) $
where $B_1, B_2, \dots$ are fixed $k\times k$ matrices.
\end{remark}

\begin{proof} \eqref{4.21} for $\ell=m$ holds by \eqref{4.10} and \eqref{4.19},
\eqref{4.20} if we note that
\begin{equation} \lb{4.22}
 \llangle P_n^*, Q_n^* \rrangle_\capL =
 \llangle Q_n, P_n \rrangle_\capR
\end{equation}
It is obvious from the definition of $\chi_\ell$ and $x_\ell$ that they lie in
the proper span. To get \eqref{4.21} for $\ell <m$, we need to know that
$\chi_\ell$ is orthogonal to the first $\ell-1$ of $1,z,z^{-1}, \dots$ and $x_\ell$
to the first $\ell-1$ of $1,z,z^{-1}, \dots$. For cases where $\chi_\ell$ or
$x_\ell$ is given by $\varphi^\capL$, this follows from  $\llangle z^k,
\varphi_\ell^\capL \rrangle =0$ for $0\leq k <\ell$ and when it is a
$(\varphi_\ell^\capR)^*$ from \eqref{4.22} which says
\[
 \llangle z^k, (\varphi_\ell^\capR)^* \rrangle_\capL =
 \llangle \varphi_\ell^\capR, z^{\ell-k} \rrangle = 0
\]
for $0\leq \ell-k <\ell$.
\end{proof}

By a (left-) module basis for $\calH_\capL$, we mean a sequence $\{f_j\}_{j=0}^\infty$
orthonormal in $\llangle \cdot \, ,\cdot\rrangle_\capL$, that is, $\llangle
f_j,f_\ell\rrangle_\calL = \delta_{j\ell}$ so that as $\{B_j\}_{j=0}^N$
runs through all $N$-tuples of $k\times k$ matrices, $\sum_{j=0}^N B_j f_j$
is a sequence of subspaces whose union is dense in $\calH_\capL$. For any
such basis, any $\eta\in\calH_\capL$ has a unique convergent expansion,
\begin{equation} \lb{4.24}
\eta =\sum_{j=0}^\infty \llangle f_j,\eta\rrangle f_j
\end{equation}

$\{\chi_j\}_{j=0}^\infty$ and $\{x_j\}_{j=0}^\infty$ are both module bases.
That means, if $\calC_{j\ell}$ is defined by
\begin{equation} \lb{4.25}
\calC_{j\ell} =\langle \chi_j, z\chi_\ell\rangle
\end{equation}
then the matrix, obtained by using the $k\times j$ blocks, $\calC_{j\ell}$,
is unitary. Moreover,
\begin{equation} \lb{4.26}
\calC_{j\ell} = \sum_m \calL_{jm} \calM_{m\ell}
\end{equation}
where
\begin{equation} \lb{4.27}
\calL_{j\ell} = \langle \chi_j, zx_\ell\rangle \qquad
\calM_{j\ell} = \langle x_j, \chi_\ell\rangle
\end{equation}

In \eqref{4.17}, set $n=2k-1$ and multiply by $z^{-k}$ to get
\begin{equation} \lb{4.28}
x_{2k} = -\alpha_n^\dag \chi_{2k} + \rho_n^\capL \chi_{2k-1}
\end{equation}
where a bottom row of $\Theta$ is clear. In this way, using \eqref{4.13},
\eqref{4.16}, \eqref{4.17}, and \eqref{4.18}, one obtains:

\begin{theorem} \lb{T4.2} With $\Theta_j (\alpha)$ given by \eqref{4.1} acting
on $\bbC^{2k}$ corresponding to $\delta_j, \delta_{j+1}$, we have
\begin{align*}
\calM &= \bdone_{1\times 1}\oplus\Theta_1(\alpha_1)\oplus\Theta_3(\alpha_3)\oplus\cdots \\
\calL &= \Theta_0(\alpha_0)\oplus\Theta_2(\alpha_2)\oplus\Theta_4(\alpha_4)\oplus\cdots
\end{align*}
\end{theorem}

The analog of \eqref{2.27} is
\begin{equation} \lb{4.29}
\calC = \begin{pmatrix}
{}& \alpha_0^\dag & \rho_0^\capL \alpha_1^\dag  &  \rho_0^\capL \rho_1 ^\capL
  & 0 & 0 & \dots & {} \\
{}& \rho_0^\capR & -\alpha_0 \alpha_1^\dag & - \alpha_0 \rho_1^\capL
  & 0 & 0 & \dots & {} \\
{}& 0 & \rho_1^\capR \alpha_2^\dag  & -\alpha_1 \alpha_2^\dag
  & \rho_2^\capL \alpha_3^\dag  & \rho_2^\capL \rho_3^\capL & \dots & {} \\
{}& 0 & \rho_1^\capR \rho_2^\capR & -\alpha_1 \rho_2^\capR
  & -\alpha_2 \alpha_3^\dag & -\alpha_2 \rho_3^\capL & \dots & {} \\
{}& 0 & 0 & 0 & \rho_3^\capR \alpha_4^\dag & -\alpha_3 \alpha_4^\dag
\end{pmatrix}
\end{equation}

We note for later purposes that for this matrix case, the GGT matrix, which we
will discuss in Section~\ref{s10}, has the form
\begin{equation} \lb{4.30}
\calG_{k\ell} = \begin{cases}
-\alpha_{k-1} \rho_k^\capL \rho_{k+1}^\capL \dots \rho_{\ell-1}^\capL \alpha_\ell^\dag
 \quad & 0 \leq k \leq \ell \\
\rho_\ell^\capR \quad & k=\ell+2 \\
0 \quad & k\geq \ell+2
\end{cases}
\end{equation}
that is,
\begin{equation} \lb{4.31}
\calG = \begin{pmatrix}
{} & \alpha_0^\dag & \rho_0^\capL \alpha_1^\dag & \rho_0^\capL \rho_1^\capL \alpha_2^\dag
  & \rho_0^\capL \rho_1^\capL \rho_2^\capL \alpha_3^\dag & {} \\
{} & \rho_0^\capR & -\alpha_0 \alpha_1^\dag & -\alpha_0 \rho_1^\capL \alpha_2^\dag
  & -\alpha_0 \rho_1^\capL \rho_2^\capL \alpha_3^\dag & {} \\
{} & 0 & \rho_1^\capR & -\alpha_1 \alpha_2^\dag
  & - \alpha_1 \rho_2^\capL \alpha_3^\dag & {} \\
{} & \dots & \dots & \dots & \dots & {}
\end{pmatrix}
\end{equation}

\section{Rank One Covariances} \lb{s5}

For selfadjoint matrices, the most elementary rank one perturbations are
diagonal, that is, $J\mapsto J+\lambda (\delta_n, \cdot\, )\delta_n$, where
$\delta_n$ is the vector with $1$ in position $n$ and $0$ elsewhere. The
impact of such a change on Jacobi parameters is trivial: $a_m\to a_m$,
$b_m\to b_m + \lambda\delta_{nm}$ (if we label vectors in the selfadjoint
case starting at $n=1$). One of our goals is to find the analog for CMV
matrices, where we will see that the impact on Verblunsky coefficients
is more subtle.

We will also address a related issue: In the spectral theory of OPUC, the
family of measures, $d\mu_\lambda$ with $\alpha_n (d\mu_\lambda)=\lambda
\alpha_n$ for a fixed $\{\alpha_n\}_{n=0}^\infty$, called an Aleksandrov
family, plays an important role analogous to a change of boundary condition
in ODE's. If $\varphi_n$ are the normalized OPUC, the GGT matrix,
\begin{equation} \lb{5.1}
\calG_{k\ell} (\{\alpha_n\}_{n=0}^\infty) = \langle\varphi_k, z\varphi_\ell\rangle
\end{equation}
has the property that $\calG (\{\lambda\alpha_n\}_{n=0}^\infty)-\calG
(\{\alpha_n\}_{n=0}^\infty)$ is rank one (see \cite[page~259]{OPUC1}). But
for the CMV basis, $\calC(\{\lambda\alpha_n\}_{n=0}^\infty) -
\calC(\{\alpha_n\}_{n=0}^\infty)$ is easily seen to be infinite rank (if the
$\alpha$'s are not mainly $0$). However, we will see here that for a suitable
$U_\lambda$ (depending on $\lambda$ but not on $\alpha$!), $U_\lambda
\calC (\{\lambda\alpha_n\}_{n=0}^\infty) U_\lambda^{-1} - \calC
(\{\alpha_n\}_{n=0}^\infty)$ is rank one.

We need to begin by figuring out what are natural rank one perturbations. The key
realization is that the proper format is multiplicative: Let $P$ be a rank one
projection and $W_\theta=e^{i\theta P} = (\bdone -P) + e^{i\theta} P$. Then $W_\theta
-\bdone = (e^{i\theta} -1)P$ is rank one, and for any $U$, $UW_\theta$ is a rank one
perturbation of $U$. It will be convenient to parametrize by $\lambda =e^{i\theta}
\in\partial\bbD$. Thus, we define
\begin{equation} \lb{5.2}
W^{(m)}(\lambda) = \bdone + (e^{i\theta} -1)(\delta_m, \cdot\, \delta_m)
\end{equation}
and given any CMV matrix $\calC$, we let
\begin{equation} \lb{5.3}
\calC^m (\lambda) = \calC W^{(m)}(\lambda)
\end{equation}
We will use $\calC^m (\lambda; \{\alpha_k\})$ where we want to make the
$\alpha$-dependence explicit. Notice that
\begin{equation} \lb{5.4}
\calC_{\ell k}^m (\lambda) =
\begin{cases}
\calC_{\ell k} \quad &\text{if } k\neq m \\
\lambda\calC_{\ell k} \quad &\text{if } k=m
\end{cases}
\end{equation}
that is, we multiply column $m$ by $\lambda$.

Part of the result we are heading towards is that
\begin{equation} \lb{5.5}
\alpha_\ell (\calC^m (\lambda)) =
\begin{cases}
\alpha_\ell (\calC) \quad & \ell <m \\
\lambda^{-1} \alpha_\ell (\calC) \quad & \ell \geq m
\end{cases}
\end{equation}
In particular, $\calC^0 (\bar\lambda)$ realizes the fact that $\calC
(\{\lambda\alpha_k\}_{k=0}^\infty)$ is unitarily equivalent to a rank one
perturbation of $\calC (\{\alpha_k\}_{k=0}^\infty)$. \eqref{5.5} for the
important case $m=0$ is due to Simon \cite[Theorem~4.2.9]{OPUC1} and
for the general case to Simon \cite{S297}. We will sketch the various proofs.

While we will eventually provide explicit unitaries that show $\calC^m
(\lambda; \{\alpha_j\}_{j=0}^\infty)$ is unitarily equivalent to $\calC$
(right side of \eqref{5.5}), we begin with a direct proof of \eqref{5.5}
in case $m=0$.

\begin{theorem}\lb{T5.1} $\calC^{m=0} (\lambda; \{\alpha_j\}_{j=0}^\infty)$
has Verblunsky coefficients $\{\lambda^{-1}\alpha_j\}_{j=0}^\infty$.
\end{theorem}

\begin{remark} If $\calM^{(\lambda)}=\lambda\bdone\oplus\Theta(\alpha_1)\oplus
\Theta(\alpha_3)\oplus\cdots$, that is, the $1$ in the upper left corner is
replaced by $\lambda$, then $\calL\calM^{(\lambda)} =\calC W^{(0)}(\lambda)$.
\end{remark}

\begin{proof}[Sketch] (See Theorems~4.2.9 and Subsection~1.4.16 of
\cite{OPUC1}.) By definition,
\begin{equation} \lb{5.6}
\calC_\lambda^{m=0} -\calC = (\lambda -1)\calC P_0
\end{equation}
where $P_0=\langle\delta_0, \cdot\,\rangle \delta_0$. Define, for $z\in\bbD$,
$F_\lambda$ and the $f_\lambda$ by
\begin{align}
F_\lambda (z) &=\langle \delta_0, (\calC^{m=0}_\lambda -z)
(C_\lambda^{m=0} +z)^{-1} \delta_0\rangle \lb{5.7} \\
&= \f{1+zf_\lambda(z)}{1-zf_\lambda (z)} \lb{5.8}
\end{align}
Using the second resolvent formula and \eqref{5.6} implies (see
Subsection~1.4.16 of \cite{OPUC1}) that
\begin{equation} \lb{5.9}
f_\lambda (z) = \lambda^{-1} f(z)
\end{equation}
The Schur algorithm and Geronimus theorem (see Chapter~3 of \cite{OPUC1}) then
imply \eqref{5.5} for $m=0$.
\end{proof}

For discussion of the movement of eigenvalues under the perturbations of
Theorem~\ref{T5.1}, see \cite{AGR,CMV2,FR} and Theorem~3.2.17 of \cite{OPUC1}.

The key to an explicit unitary equivalence is the following. Let
\begin{equation} \lb{5.10}
\nu(\lambda) = \begin{pmatrix}
1 & 0 \\ 0 & \lambda \end{pmatrix}
\qquad
\ti\nu (\lambda) = \begin{pmatrix}
\lambda & 0 \\ 0 & 1 \end{pmatrix}
\end{equation}
Then, by a simple calculation,
\begin{align}
\nu(\lambda)\Theta (\lambda^{-1}\alpha) \nu(\lambda)
  &= \lambda \Theta (\alpha) \lb{5.11} \\
\ti\nu(\lambda)^{-1}\Theta (\lambda^{-1} \alpha) \ti\nu (\lambda)^{-1}
  &= \lambda^{-1} \Theta (\alpha) \lb{5.12}
\end{align}
Note that \eqref{5.11} does not use $\nu(\lambda)$ and $\nu (\lambda)^{-1}$
but $\nu (\lambda)$ in both places. Similarly, \eqref{5.12} has $\ti\nu
(\lambda)^{-1}$ in both places. In the full calculation, one does not use
\[
U\!\calL\calM U^{-1} = (U\!\calL U^{-1})(U\!\calM U^{-1})
\]
but rather
\begin{equation} \lb{5.13}
U\!\calL\calM U^{-1} = (U\!\calL U)(U^{-1}\calM U^{-1})
\end{equation}

We need a notation for diagonal matrices. $D (1^{2k} (1\lambda)^\infty)$ indicates
the diagonal matrix with entries $1$ $2k$ times, then alternating $1$'s and
$\lambda$'s. Thus,
\begin{equation} \lb{5.14}
W^{(m)}(\lambda) = D(1^m \lambda 1^\infty)
\end{equation}
Using \eqref{5.14}, \eqref{5.13}, \eqref{5.11}, and \eqref{5.12}, a direct
calculation (see Section~5 of \cite{S297}) shows:

\begin{theorem} \lb{T5.2} For $n=0,1,2,\dots$, define
\begin{gather}
U_{2k-1} = D(1^{2k} (1\lambda)^\infty) \lb{5.15a} \\
U_{2k} = D(\lambda^{2k} (1\lambda)^\infty) \lb{5.15b} \\
T_{n,\lambda} (\{\alpha_j\}_{j=0}^\infty) = \beta_j \lb{5.16}
\end{gather}
where
\begin{equation} \lb{5.17}
\beta_j = \begin{cases}
\alpha_j \quad & j<n \\
\lambda\alpha_j \quad & j\geq n
\end{cases}
\end{equation}
Then
\begin{equation} \lb{5.18}
U_n\calC (T_{n,\lambda^{-1}} (\{\alpha_j\}_{j=0}^\infty)) U_n^{-1} =
\calC (\{\alpha_j\}_{j=0}^\infty) W^{(n)} (\lambda)
\end{equation}
In particular, \eqref{5.5} holds.
\end{theorem}

\begin{remarks} 1. It is important that $\delta_0$ is an eigenvector of $U_n$
since Verblunsky coefficients involve a unitary and a cyclic vector. \eqref{5.18}
also shows that $\calC W^{(n)}(\lambda)$ has $\delta_0$ as a cyclic vector.

\smallskip
2. One can also ask about Verblunsky coefficients of $W^{(n)}(\lambda) \calC
(\{\alpha_j\}_{j=0}^\infty)$. Since Verblunsky coefficients are invariant under
unitaries that have $\delta_0$ as an eigenvector and
\[
W^{(n)}\calC = W^{(n)}\calC W^{(n)} (W^{(n)})^{-1}
\]
the Verblunsky coefficients of $\calC W^{(n)}$ and $W^{(n)}\calC$ are the same.
\end{remarks}

\eqref{5.11}, \eqref{5.12}, and \eqref{5.13} imply a result about extended
CMV matrices. For $\lambda\in\partial\bbD$, let $\wti W(\lambda)$ be the two-sided
diagonal matrix with $d_{2j}=1$, $d_{2j+1}=\lambda$. Then

\begin{theorem}\lb{T5.3} Let $\lambda\in\partial\bbD$. Then $\wti W(\lambda)
\calE(\{\alpha_n\}_{n=-\infty}^\infty)\wti W(\lambda)^{-1} =\calE
(\{\lambda \alpha_n\}_{n=-\infty}^\infty)$.
\end{theorem}

\begin{remark} In particular, spectral properties of $\calE
(\{\alpha_n\}_{n=-\infty}^\infty)$ and $\calE (\{\lambda\alpha_n\}_{n=-\infty}^\infty)$
are identical and $\alpha_n \to \lambda\alpha_n$ preserves isospectral tori in
the periodic case.
\end{remark}

\section{Resolvents of CMV Matrices} \lb{s6}

In this section, we will present formulae for the resolvent of $\calC$ analogous to
the Green's function formula for Jacobi matrices (see Theorem~4.4.3 of \cite{OPUC1}).
These formulae appeared first in Section~4.4 of \cite{OPUC1}. Similar formulae for
GGT matrices appeared earlier in Geronimo--Teplyaev \cite{GTep} (see also
\cite{GJo2,GJo1}).

Clearly, we need an analog of Jost solutions. For OPUC, these were found by
Golinskii--Nevai \cite{GN} who proved

\begin{theorem} \lb{T6.1} Fix $z\in\bbD$. Let $\varphi_n$ be the normalized OPUC
for a probability measure $d\mu$ on $\partial\bbD$, and $\psi_n$ the normalized
OPUC for Verblunsky coefficients $-\alpha_n (d\mu)$ {\rm{(}}so-called second
kind polynomials{\rm{)}}. Then
\begin{equation} \lb{6.1}
\sum_{n=0}^\infty \, \abs{\psi_n(z) + F(z)\varphi_n(z)}^2 +
\abs{\psi_n^*(z) - F(z) \varphi_n^*(z)}^2 <\infty
\end{equation}
where $F$ is the Carath\'eodory function:
\begin{equation} \lb{6.2}
F(z)\equiv \int \f{e^{i\theta} + z}{e^{i\theta} - z}\, d\mu(\theta)
\end{equation}
\end{theorem}

\begin{remarks} 1. This is an analog of Weyl's formula; see (1.2.53) of \cite{OPUC1}.

\smallskip
2. See \cite{GN} or \cite[Section~3.2]{OPUC1} for a proof.
\end{remarks}

With this in mind, we define
\begin{align}
y_n &= \begin{cases} z^{-\ell} \psi_{2\ell} \quad & n=2\ell \\
-z^{-\ell} \psi_{2\ell-1}^* \quad & n=2\ell-1 \end{cases} \lb{6.3} \\
\notag \\
\Upsilon_n &= \begin{cases} -z^{-\ell} \psi_{2\ell}^* \quad & n=2\ell \\
z^{-\ell +1} \psi_{2\ell -1} \quad & n=2\ell-1 \end{cases} \lb{6.4} \\
\notag \\
p_n &= y_n + F(z) x_n \lb{6.5} \\
\pi_n & = \Upsilon_n + F(z) \chi_n \lb{6.6}
\end{align}
Then Theorem~4.4.1 of \cite{OPUC1} says:

\begin{theorem} \lb{T6.2} We have that for $z\in\bbD$,
\begin{equation} \lb{6.6a}
[(\calC-z)^{-1}]_{k\ell} =
\begin{cases} (2z)^{-1}\chi_\ell(z) p_k (z) \quad & k>\ell \text{ or } k=\ell =2n-1 \\
(2z)^{-1} \pi_\ell (z) x_k (z) \quad & \ell >k \text{ or } k=\ell=2n \end{cases}
\end{equation}
\end{theorem}

As a special case, since $\delta_n =\chi_n(\calC)\delta_0$ and $\abs{\chi_n (e^{i\theta})}
= \abs{\varphi_n (e^{i\theta})}$, we obtain from a spectral representation
\begin{equation} \lb{6.7}
\int \f{\abs{\varphi_n (e^{i\theta})}^2}{e^{i\theta} -z}\, d\mu(\theta) =
(2z^{n+1})^{-1} \varphi_n(z) [-\psi_n^* (z) + F(z) \varphi_n^*(z)]
\end{equation}

As shown in remarks to Theorem~9.2.4 of \cite{OPUC2}, this is equivalent to a formula
of Khrushchev \cite{Kh2000} for $\int \f{e^{i\theta}+z}{e^{i\theta}-z} \abs{\varphi_n
(e^{i\theta})}^2\,d\mu (\theta)$. For an application of Theorem~\ref{T6.2},
see Stoiciu \cite{StoiJAT}.

\section{$\calI^p$ Perturbations} \lb{s7}

In this section, we give some elementary estimates of Golinskii--Simon \cite{GSunpub}
on the $\calI^p$ norm of $\calC (\{\alpha_n\}_{n=0}^\infty)- \calC(\{\beta_n\}_{n=0}^\infty)$.
(For definition and background on Schatten $p$-classes, see Gohberg--Krein \cite{GK}
and Simon \cite{STI}.)

In this section, we will use these estimates to write the Szeg\H{o} function in terms
of Fredholm determinants of the CMV matrices and to discuss scattering theory.
Further applications appear in Section~\ref{s9}.

A diagonal matrix, $A$, has $\calI_p$ norm $(\sum_j \abs{a_{jj}}^p)^{1/p})$ and $\|A\|_p$
is invariant under multiplication by a unitary. So if $A$ has only a nonvanishing $k$th
principal diagonal, $A$ has $\calI_p$ norm $(\sum_j \abs{a_{j\, j+k}}^p)^{1/p}$. Since
$(a^{1/p} + b^{1/p} + c^{1/p}) \leq (a+b+c)^{1/p} 3^{1-1/p}$ (by H\"older's inequality),
we see for tridiagonal matrices that
\begin{equation} \lb{7.1}
\|A-B\|_p \leq 3^{1-1/p} \biggl(\, \sum_{i,j} \abs{a_{ij}-b_{ij}}^p\biggr)^{1/p}
\end{equation}
This lets us estimate $\|\calL(\{\alpha_n\}_{n=0}^\infty) -
\calL(\{\beta_n\}_{n=0}^\infty)\|_p$, and similarly for $\calM$. Using unitarity of $\calL$
and $\calM$, $\|\calL\calM-\calL'\calM'\|_p\leq \|(\calL-\calL')\calM\|_p +
\|\calL' (\calM-\calM')\|_p \leq \|\calL-\calL'\|_p + \|\calM-\calM'\|_p$. So using
$(a^{1/p} + b^{1/p})\leq 2^{1-1/p} (a+b)^{1/p}$, we find

\begin{theorem}[$=$ Theorem~4.3.2 of \cite{OPUC1}]\lb{T7.1} Let $\{\alpha_n\}_{n=0}^\infty$,
$\{\beta_n\}_{n=0}^\infty$ be two sequences in $\ol{\bbD}^\infty$ and let $\rho_n =
(1-\abs{\alpha_n}^2)^{1/2}$, $\sigma_n = (1-\abs{\beta_n}^2)^{1/2}$. Then
\begin{equation} \lb{7.2}
\|\calC(\{\alpha_n\}_{n=0}^\infty) - \calC(\{\beta_n\}_{n=0}^\infty)\|_p \leq 6^{1-1/p}
\biggl(\, \sum_{n=0}^\infty\, \abs{\alpha_n -\beta_n}^p +
\abs{\rho_n -\sigma_n}^p\biggr)^{1/p}
\end{equation}
\end{theorem}

\begin{remark} \cite{OPUC1} has the constants $2$ for $1\leq p\leq 2$ and $2\cdot 3^{1-1/2p}$
for $2\leq p\leq\infty$, but the proof there actually shows $2^{1-1/p}$ and $2^{1-1/p}\,
3^{1-1/2p}$. This improves the constant in \eqref{7.2}.
\end{remark}

To rephrase in terms of $\abs{\alpha_n -\beta_n}$ only, we first note that
\[
\sup_{\abs{z}\leq R}\, \biggl|\, \f{d}{dz}\, (1-\abs{z}^2)^{1/2} \biggr|
\leq (1-R^2)^{-1/2}
\]
and $\abs{\abs{z}-\abs{w}} \leq \abs{z-w}$ to see that
\begin{equation} \lb{7.3}
\sup_{\abs{z},\abs{w}\leq R}\, \abs{(1-\abs{z}^2)^{1/2} - (1-\abs{w}^2)^{1/2}}
\leq (1-R^2)^{-1/2} \abs{z-w}
\end{equation}
We need to note that $\abs{\sqrt{a} -\sqrt{b}\,} \leq \sqrt{\abs{a-b}}$ and
$\abs{\abs{\alpha}^2 -\abs{\beta}^2}\leq 2 \abs{\alpha-\beta}$, so
\begin{equation} \lb{7.4}
\abs{(1-\abs{z}^2)^{1/2} - (1-\abs{w}^2)^{1/2}} \leq \sqrt{2}\,
\abs{z-w}^{1/2}
\end{equation}
Thus,

\begin{theorem}\lb{T7.2} Let $\{\alpha_n\}_{n=0}^\infty$ and $\{\beta_n\}_{n=0}^\infty$
be two sequences in $\ol{\bbD}^\infty$ and let $\rho_n =(1-\abs{\alpha_n}^2)^{1/2}$,
$\sigma_n = (1-\abs{\beta_n}^2)^{1/2}$. Then
\begin{SL}
\item[{\rm{(a)}}] If $\sup_n \abs{\alpha_n} \leq R<1$ and $\sup_n \abs{\beta_n}\leq R$,
then
\begin{equation} \lb{7.5}
\|\calC(\{\alpha_n\}_{n=0}^\infty) - \calC(\{\beta_n\}_{n=0}^\infty)\|_p \leq 6^{1-1/p}
[1+(1-R^2)^{-p/2}]^{1/p} \biggl(\, \sum_{n=0}^\infty\,
\abs{\alpha_n -\beta_n}^p \biggr)^{1/p}
\end{equation}

\item[{\rm{(b)}}] In general, for $1\leq p\leq \infty$,
\begin{equation} \lb{7.6}
\|\calC(\{\alpha_n\}_{n=0}^\infty) - \calC(\{\beta_n\}_{n=0}^\infty)\|_p \leq 6^{1-1/p}
\biggl(\, \sum_{n=0}^\infty\, \abs{\alpha_n -\beta_n}^p + 2^{p/2}
\abs{\alpha_n -\beta_n}^{p/2} \biggr)^{1/p}
\end{equation}
\end{SL}
\end{theorem}

One thing made possible by CMV matrices is scattering theory because all CMV
matrices act on the same space $(\ell^2 (\{0,1,2,\dots\}))$. This is an important
tool made possible by the CMV matrix. Thus,

\begin{theorem} \lb{T7.3} Suppose $\sup_n \abs{\alpha_n}\leq R \leq 1$,
$\sup_n \abs{\beta_n}\leq R <1$, and
\begin{equation} \lb{7.6x}
\sum_{n=0}^\infty \, \abs{\alpha_n-\beta_n} <\infty
\end{equation}
Then, if $P_\ac (\cdot)$ is the projection onto the absolutely continuous
subspace of an operator and $\calC=\calC(\{\alpha_n\}_{n=0}^\infty)$, $\wti\calC
=\calC (\{\beta_n\}_{n=0}^\infty)$, then
\[
\lim_{n\to\pm\infty} \, \calC^n \wti\calC^{-n} P_\ac (\wti\calC)
\]
exists and is a partial isometry with range $P_\ac (\calC)$. In particular, $\calC
\restriction  P_\ac (\calC)$ and $\wti\calC\restriction P_\ac (\wti\calC)$ are
unitarily equivalent.
\end{theorem}

\begin{remarks} 1. This follows from the fact that $\calC-\wti\calC$ is trace
class and from the Kato--Birman theorem \cite{RS3}.

\smallskip
2. If $\{\alpha_n\}_{n=0}^\infty$ corresponds to
\begin{equation} \lb{7.6a}
d\mu = f(\theta)\, \f{d\theta}{2\pi} + d\mu_\delta
\end{equation}
and $\{\beta_n\}_{n=0}^\infty$ corresponds to
\begin{equation} \lb{7.6b}
d\nu = g(\theta)\, \f{d\theta}{2\pi} +  d\nu_\delta
\end{equation}
then this theorem implies that up to sets of $d\theta$-measure $0$,
\begin{equation} \lb{7.7}
\{\theta\mid f(\theta)\neq 0\} = \{\theta\mid g(\theta)\neq 0\}
\end{equation}
(also see Theorem~\ref{T9.3}).

\smallskip
3. For the case $\beta_n \equiv 0$, this holds if only  $\sum_{n=0}^\infty
\abs{\alpha_n}^2 <\infty$; see \cite[Section~10.7]{OPUC2}.
\end{remarks}

Finally, following Simon \cite[Section~4.2]{OPUC1}, we want to state the connection
of $\calC$ to the Szeg\H{o} function, defined for $\abs{z}<1$ by
\begin{equation} \lb{7.8}
D(z) = \lim_{n\to\infty}\, \varphi_n^* (z)^{-1}
\end{equation}
which exists and is nonzero if (and only if)
\begin{equation} \lb{7.9}
\sum_{j=0}^\infty\, \abs{\alpha_j}^2 <\infty
\end{equation}
(see Section~2.4 of \cite{OPUC1}). We will let $\calC_0$ be the free CMV matrix
corresponding to $d\mu =\f{d\theta}{2\pi}$; equivalently, $\alpha_n\equiv 0$.

\begin{theorem}\lb{T7.4} Suppose
\begin{equation} \lb{7.10}
\sum_{n=0}^\infty\, \abs{\alpha_n} <\infty
\end{equation}
Then $\calC-\calC_0$ is trace class and
\begin{equation} \lb{7.11}
D(z)^{-1} D(0) = \det \biggl( \f{1-z\bar\calC}{1-z\bar\calC_0}\biggr)
\end{equation}
If \eqref{7.9} holds, then $\calC-\calC_0$ is Hilbert--Schmidt, and
\begin{equation} \lb{7.12}
D(z)^{-1} D(0) = \det_2 \biggl( \f{1-z\bar\calC}{1-z\bar\calC_0}\biggr) e^{-zw_1}
\end{equation}
where
\begin{equation} \lb{7.13}
w_1 = \alpha_0 -\sum_{j=1}^\infty \alpha_j \bar\alpha_{j-1}
\end{equation}
\end{theorem}

\begin{remarks} 1. Alas, (4.2.53) of \cite{OPUC1} has a sign error: it is
$e^{-zw_1}$ as we have here, not $e^{zw_1}$ as appears there!

\smallskip
2. By $\det(\f{1-z\bar\calC}{1-z\bar\calC_0})$, we mean $\det((1-z\bar\calC)
(1-z\bar\calC_0)^{-1})$. Since
\[
(1-z\bar\calC)(1-z\bar\calC_0)^{-1} = 1-z (\bar\calC - \bar\calC_0)
(1-z\bar\calC_0)^{-1}
\]
we see that this is $1+$ trace class (resp., Hilbert--Schmidt) if
$\calC-\calC_0$ is trace class (resp., Hilbert--Schmidt).

\smallskip
3. For a proof, see Section~4.2 (Theorem~4.2.14) of \cite{OPUC1}.

\smallskip
4. $\bar\calC$ is the complex conjugate of $\calC$, that is, $(\bar\calC)_{ij}
=\ol{(\calC_{ij})}$.

\smallskip
3. $\det(\cdot)$ is defined on operators of the form $1+A$ with $A$ trace class,
and then $\det_2$ on $1+A$ with $A$ Hilbert--Schmidt by
\begin{equation} \lb{7.14}
\det (1+A) = \det_2 ((1+A)e^{-A})
\end{equation}
When $A$ is trace class,
\begin{equation} \lb{7.15}
\det (1+A) = \det_2 (1+A) e^{\tr (A)}
\end{equation}
If \eqref{7.10} holds, $-zw_1 =\tr ((1-z\bar\calC)/(1-z\bar\calC_0))$ and
\eqref{7.11}/\eqref{7.12} are consistent by \eqref{7.15}. See \cite{GK} or
\cite{STI} for a discussion of $\det(\cdot)$ and $\det_2 (\cdot)$.

\smallskip
5. The connection for one-dimensional Schr\"odinger operators of the Jost
function and Fredholm determinants goes back to Jost--Pais \cite{JP}. For Jacobi
matrices, under the name ``perturbation determinant," they were used by
Killip--Simon \cite{S281}.
\end{remarks}

\section{Essential Spectra} \lb{s8}

The {\it discrete spectrum\/} of an operator is the set of isolated points of finite
multiplicity. The complement of the discrete spectrum in the spectrum is called the
{\it essential spectrum}. Since a CMV matrix has a cyclic vector, the essential
spectrum is just the set of nonisolated points in the support of the spectral
measure, $d\mu$, often called the {\it derived set\/} of $\supp (d\mu)$.
Last--Simon \cite{LSjdam} have a general result for the essential spectrum of
a CMV matrix $\calC (\{\alpha_n\}_{n=0}^\infty)$.

\smallskip
\noindent{\bf Definition.} A {\it right limit\/} of $\{\alpha_n\}_{n=0}^\infty$ is any
two-sided sequence $\{\beta_n\}_{n=-\infty}^\infty$ in $\ol{\bbD}^\bbZ$ for which
there exists $n_\ell\to\infty$ so $\lim_{\ell\to\infty} \alpha_{n_\ell +j} =
\beta_j$ for each $j\in\bbZ$. $\calR(\{\alpha_n\}_{n=0}^\infty)$ is the set
of right limits of $\{\alpha_n\}_{n=0}^\infty$.

Since $\ol{\bbD}^\bbZ$ is compact, $\calR$ is nonempty. Indeed, if $\ti\beta_0$
is any limit point of $\alpha_n$, there is a right limit with $\beta_0 =
\ti\beta_0$.

\begin{theorem}[Last--Simon \cite{LSjdam}] \lb{T8.1} For any $\{\alpha_n\}_{n=0}^\infty
\in\ol{\bbD}^\infty$, we have
\begin{equation} \lb{8.1}
\sigma_\ess (\calC(\{\alpha_n\}_{n=0}^\infty)) =
\ol{\bigcup_{\beta\in\calR (\{\alpha_n\}_{n=0}^\infty)}\, \sigma (\calE(\{\beta_n\}_{n=0}^\infty))}
\end{equation}
\end{theorem}

\begin{remarks} 1. The proof \cite{LSjdam} uses a Weyl trial sequence argument.
The key is that because $\calC$ has finite width, for any $\lambda_0\in\sigma_\ess(\calC)$
and $\veps$, there exist $L$, $n_j\to\infty$ and $\varphi_j$ supported in $(n_j-L, n_j+L)$
with $\|\varphi_j\|=1$ and
\begin{equation} \lb{8.2}
\limsup_{j\to\infty}\, \|(\calC-\lambda_0) \varphi_j\|\leq \veps
\end{equation}

\smallskip
2. Right limits of Verblunsky coefficients were considered earlier by Golinskii--Nevai
\cite{GN}, motivated by earlier work on Schr\"odinger operators by Last--Simon \cite{LastS}.
This work was in the context of a.c.\ spectrum (see Theorem~10.9.11(ii) of \cite{OPUC2}).

\smallskip
3. \cite{LSjdam} used the same methods to study Jacobi and Schr\"odinger operators. Earlier
results of the form \eqref{7.1} for Schr\"odinger operators (but not for CMV matrices)
are due to Georgescu--Iftimovici \cite{GI1}, M\u antoiu \cite{M1}, and
Rabinovich \cite{Rab05}. These rely on what I regard as elaborate machines
(connected with $C^*$-algebras or with Fredholm operators) although, no doubt,
their authors regard them as very natural.
\end{remarks}

One can combine this with Theorem~\ref{T5.3} to obtain

\begin{theorem} \lb{T8.2} Let $\{\alpha_j\}_{j=0}^\infty$ and $\{\beta_j\}_{j=0}^\infty$
be two sequences of Verblunsky coefficients. Suppose there exist $\lambda_j\in\partial\bbD$
so that
\begin{alignat}{2}
\text{\rm{(i)}}& \qquad &\beta_j\lambda_j -\alpha_j &\to 0  \lb{8.3} \\
\text{\rm{(ii)}} &  \qquad & \lambda_{j+1} \bar\lambda_j &\to 1  \lb{8.4}
\end{alignat}
Then
\begin{equation} \lb{8.5}
\sigma_\ess (\calC (\{\alpha_j\}_{j=0}^\infty)) =
\sigma_\ess (\calC (\{\beta_j\}_{j=0}^\infty ))
\end{equation}
\end{theorem}

\begin{proof} Let $\{\gamma_j\}_{j=-\infty}^\infty$ be a right limit of $\{\beta_j\}_{j=0}^\infty$.
By passing to a subsequence, we can suppose $\lambda_{n_j}\to\lambda_\infty$ and $\beta_{n_j+k}
\to \gamma_k$. Since $\lambda_{n_j+k}\lambda_{n_j}^{-1}\to 1$, we see that $\alpha_{n_j+k}
\to \lambda_\infty\gamma_k$. By Theorem~\ref{T5.3}, $\sigma (\calE(\{\gamma_k\}_{k=-\infty}^\infty))
=\sigma (\calE (\{\lambda_\infty\gamma_k\}_{k=-\infty}^\infty))$. It follows (using symmetry)
that \eqref{8.5} holds.
\end{proof}

\begin{remark} This proof is from \cite{LSjdam}, but the result appears earlier as Theorem~4.3.8
in \cite{OPUC1}, motivated by a special case of Barrios--L\'opez \cite{BRLL}.
\end{remark}

\begin{example}\lb{E8.3} (This is due to Golinskii \cite{Gol2000}; the method of proof is due
to \cite{LSjdam}. See the discussion in \cite{LSjdam} for earlier related results.) Suppose
$\abs{\alpha_n}\to 1$ as $n\to\infty$. Then $\sigma_\ess (\calC(\{\alpha_n\})_{n=0}^\infty)$
is the set of limit points of $\{-\bar\alpha_{j+1}\alpha_j\}_{j=0}^\infty$. For any limit
point has $\calE(\{\beta_j\}_{j=0}^\infty)$ diagonal (since $(1-\abs{\beta_j}^2)^{1/2}\equiv 0$)
with diagonal values $-\bar\beta_{j+1}\beta_j$, and by compactness, any limit point of
$-\bar\alpha_{j+1}\alpha_j$ occurs as some $-\bar\beta_1 \beta_0$. In particular, this (plus
an extra argument) implies $\sigma_\ess (\calE)=\{\lambda_0\}$ if and only if $\abs{\alpha_n}
\to 1$ and $\bar\alpha_{n+1}\alpha_n\to -\lambda_0$. See \cite{Gol2000} and
\cite{LSjdam} for a discussion of when $\sigma_\ess (\calC)$ is a finite set. \qed
\end{example}

It is well known (see Example~1.6.12 of \cite{OPUC1} and Example~11.1.4 of \cite{OPUC2})
that if $\alpha_n\equiv a\in\bbD$, then $\sigma_\ess (\calC)=\Delta_{\abs{a}}=
\{z\in \partial\bbD\mid \abs{\arg z}\geq 2 \arcsin (\abs{a})\}$, which increases
as $\abs{a}$ decreases. It follows:

\begin{example}\lb{E8.4}  (Theorem~7.8 of \cite{LSjdam}; one direction was proven in
\cite{CMVppt}, which motivated Theorem~7.8 of \cite{LSjdam}.) Suppose
\begin{equation} \lb{8.6}
\f{\alpha_{j+1}}{\alpha_j}\to 1 \qquad \liminf \abs{\alpha_j} =a
\end{equation}
Then
\begin{equation} \lb{8.7}
\sigma_\ess (\calC(\{\alpha_n\}_{n=0}^\infty)) =\Delta_a
\end{equation}
For $\f{\alpha_{j+1}}{\alpha_j} \to 1$ implies that each limit is of the form
$\beta_j \equiv b$ for some $b\in\bbD$, so
\[
\sigma_\ess (\calC) =\bigcup_{b=\text{limits of }\alpha_j }\, \Delta_{\abs{b}}
= \Delta_a
\]
since $\Delta_{\abs{b}}\subseteq \Delta_a$ if $\abs{b}\geq a$.  \qed
\end{example}

\section{Spectral Consequences} \lb{s9}

Section~4.3 of \cite{OPUC1} describes joint work of Golinskii--Simon \cite{GSunpub}
that uses CMV matrices to obtain spectral results that relate properties of
$\{\alpha_n\}_{n=0}^\infty$ to the associated measures. Here, in brief, are some
of their main results:

\begin{theorem}[$\equiv$ Theorem~4.3.5 of \cite{OPUC1}; subsumed in Theorem~\ref{T8.2}]
\lb{T9.1} If $\abs{\alpha_n-\beta_n} \to 0$, then $\sigma_\ess (\calC(\{\alpha_n\}_{n=0}^\infty))
=\sigma_\ess (\calC(\{\beta_n\}_{n=0}^\infty))$.
\end{theorem}

\begin{remark} Of course, Theorem~\ref{T9.1} also follows from Theorem~\ref{T8.1}.
\end{remark}

\begin{proof} By \eqref{7.6} and a limiting argument, $\calC(\{\alpha_n\}_{n=0}^\infty)
- \calC (\{\beta_n\}_{n=0}^\infty)$ is compact. The result follows from Weyl's theorem
on the invariance of essential spectrum under compact perturbation.
\end{proof}

\begin{theorem}[$\equiv$ Theorem~4.3.4 of \cite{OPUC1}]\lb{T9.2} If $\limsup
\abs{\alpha_n(d\mu)}=1$, then $d\mu$ is purely singular.
\end{theorem}

\begin{remark} This result is called Rakhmanov's lemma, after \cite{Rakh83}. The proof
is motivated by earlier results for Jacobi matrices of Dombrowski \cite{Dom78} and
Simon--Spencer \cite{SimSp}.
\end{remark}

\begin{proof} Let $\widehat\alpha_n$ be defined by
\[
\widehat\alpha_n = \begin{cases}
1 \quad &\text{if } \alpha_n =0\\
\f{\alpha_n}{\abs{\alpha_n}} \quad &\text{if } \alpha_n\neq 0
\end{cases}
\]
Since $\limsup\abs{\alpha_n}=1$, we can find a sequence $n_j\to\infty$, so
\[
\sum_{j=0}^\infty \, \abs{\alpha_{n_j}-\widehat\alpha_{n_j}}^{1/2} < \infty
\]
Let
\[
\beta_n =\begin{cases}
\widehat\alpha_n \quad &\text{if $n=n_j$ for some $j$} \\
\alpha_n \quad &\text{otherwise}
\end{cases}
\]
Then $\calC(\{\beta_n\}_{n=0}^\infty) - \calC(\{\alpha_n\}_{n=0}^\infty)$ is trace
class by \eqref{7.6}. By the Kato--Birman theorem \cite{RS3},
\[
\sigma_\ac (\calC(\{\alpha_n\}_{n=0}^\infty)) =
\sigma_\ac (\calC (\{\beta_n\}_{n=0}^\infty))
\]
Since $\abs{\widehat\alpha_n}=1$, $\calC(\{\beta_n\}_{n=0}^\infty)$ is a direct
sum of finite matrices of size $n_{j+1}-n_j$, and so it has no a.c.\ spectrum.
\end{proof}

\begin{theorem}[$\equiv$ Theorem~4.3.6 of \cite{OPUC1}] \lb{T9.3} If
$\{\alpha_n\}_{n=0}^\infty$ and $\{\beta_n\}_{n=0}^\infty$ are the Verblunsky
coefficients of $d\mu$ and $d\nu$ given by \eqref{7.6a} and \eqref{7.6b}, and
\eqref{7.6x} holds, then \eqref{7.7} holds.
\end{theorem}

\begin{proof} If $\limsup \abs{\alpha_n} <1$, then $\limsup \abs{\beta_n} <1$,
and by Theorem~\ref{T7.3}, \eqref{7.7} holds. If $\limsup \abs{\alpha_n}
=\limsup \abs{\beta_n}=1$, then $\sigma_\ac (\calC(\{\alpha_n\}_{n=0}^\infty))
= \sigma_\ac (\calC(\{\beta_n\}_{n=0}^\infty))=\emptyset$ by Theorem~\ref{T9.2}.
\end{proof}

\section{The AGR Factorization of GGT Matrices} \lb{s10}

This section is primarily preparatory for the next and discusses GGT matrix
representations (for Geronimus \cite{Ger44}, Gragg \cite{Gragg82}, and
Teplyaev \cite{Tep91}) associated to a measure on $\partial\bbD$:
\begin{equation} \lb{10.1}
\calG_{k\ell} (\{\alpha_n\}_{n=0}^M) = \langle \varphi_k, z\varphi_\ell\rangle
\end{equation}
and discussed in Section~4.1 of \cite{OPUC1}. If $\mu$ is nontrivial, $M$ in
\eqref{10.1} is $\infty$ and $\alpha_n\in\bbD$ for all $n$. If $\mu$ is supported
on exactly $N$ points, $M=N-1$ and $\alpha_0, \dots, \alpha_{N-2}\in\bbD$,
$\alpha_{N-1}\in\partial\bbD$. There is an explicit calculation (see Proposition~1.5.9
of \cite{OPUC1}):
\begin{equation} \lb{10.2}
\calG_{k\ell} = \begin{cases}
-\bar\alpha_\ell \alpha_{k-1} \prod_{m=k}^{\ell-1} \rho_m
\quad & 0\leq k\leq \ell \\
\rho_\ell \quad & k=\ell+1 \\
0 \quad & k\geq \ell+2
\end{cases}
\end{equation}

We present a remarkable factorization of GGT matrices due to Ammar, Gragg,
and Reichel \cite{AGR}, use it to provide a result about cosets in
$\bbU(N)/\bbU(N-1)$, and then present an alternate proof of Theorems~\ref{T9.2}
and \ref{T9.3} using GGT rather than CMV matrices. For the special case of
orthogonal matrices (all $\alpha_j\in (-1,1)$), AGR found this factorization
earlier \cite{Athens}.

We defined $\Theta_j(\alpha)$ before \eqref{3.13} as a $2\times 2$ matrix acting
on the span of $\delta_j, \delta_{j+1}$. We define $\wti\Theta (\alpha_j)$ to
be this matrix viewed as an operator on $\bbC^N$ by $\bdone_j\oplus\Theta_j
(\alpha)\oplus\bdone_{N-j-2}$. $\wti{\wti\Theta}_{N-1}(\alpha)$ is the matrix
$\bdone_{N-1}\oplus\bar\alpha$.

\begin{theorem}[AGR factorization]\lb{T10.1} For any finite $N\times N$ GGT matrix,
\begin{equation} \lb{10.3}
\calG (\{\alpha\}_{n=0}^{N-1}) = \wti\Theta_0 (\alpha_0) \dots \wti\Theta_{N-2}
(\alpha_{N-2}) \wti{\wti \Theta}_{N-1} (\alpha_{N-1})
\end{equation}
For $N=\infty$,
\begin{equation} \lb{10.4}
\calG (\{\alpha_n\}_{n=0}^\infty) = \slim_{M\to\infty}\, \wti\Theta_0 (\alpha_0)
\dots\wti\Theta_M (\alpha_M)
\end{equation}
\end{theorem}

\begin{remarks} 1. \eqref{10.4} follows from \eqref{10.3} by a simple limiting
argument. We will only prove \eqref{10.3} below.

\smallskip
2. We will give three proofs which illustrate slightly different aspects of the
formula.

\smallskip
3. As explained in Killip--Nenciu \cite{KN}, the Householder algorithm lets
one write any unitary as a product of $N-1$ reflections; in many ways, the
representation \eqref{10.4} is more useful.
\end{remarks}

\begin{proof}[First Proof] By a direct calculation using \eqref{10.2},
\begin{equation} \lb{10.5x}
\calG (\{\alpha_n\}_{n=0}^{N-1}) =\Theta_0(\alpha_0) [\bdone_{1\times 1} \oplus
\calG(\{\alpha_{n+1}\}_{n=0}^{N-2})]
\end{equation}
\eqref{10.3} follows by induction.
\end{proof}

\begin{proof}[Second Proof {\rm{(}}that of {\rm{AGR \cite{AGR})}}] We will prove
first that any unitary (upper) Hessenberg matrix $H$ (i.e., $H_{k\ell}=0$ if
$k\geq \ell+1$) with positive subdiagonal (i.e., $H_{\ell+1,\ell}>0$ for all $\ell$)
has the form \eqref{10.3} for suitable $\alpha_0, \alpha_1, \dots, \alpha_{N-2}\in
\bbD$ and $\alpha_{N-1}\in\partial\bbD$. The first column of $H$ has the form
$(\bar\alpha_0,\rho_0, 0, \dots, 0)^t$ for some $\alpha_0\in\bbD$. Then
$\Theta_0 (\alpha_0)^{-1}H$ is of the form $\bdone_{1\times 1}\oplus H^{(1)}$,
where $H^{(1)}$ is a unitary $(N-1)\times (N-1)$ Hessenberg matrix with positive
subdiagonal. By induction, $H$ has the form \eqref{10.3}. One proves that
$\{\alpha_n\}_{n=0}^{N-1}$ are the Verblunsky coefficients of the GGT matrix,
either by using \eqref{10.2} or by deriving recursion relations.
\end{proof}

For the third proof, we need a lemma that is an expression of Szeg\H{o} recursion.

\begin{lemma} \lb{L10.2} We have that
\begin{align}
\langle \varphi_j^*, z\varphi_j\rangle &= \bar\alpha_j \lb{10.4a} \\
\langle \varphi_{j+1}, z\varphi_j\rangle &=  \rho_j \lb{10.4b} \\
\langle \varphi_{j+1}, \varphi_{j+1}^* \rangle &= -\bar\alpha_j \lb{10.4c} \\
\langle \varphi_j^*, \varphi_{j+1}^*\rangle &=  \rho_j \lb{10.4d}
\end{align}
\end{lemma}

\begin{remark} This says that a certain change of basis on a two-dimensional
space is $\Theta (\alpha_j)$.
\end{remark}

\begin{proof} $\varphi_{j+1}\perp\varphi_j^*$ since $\deg (\varphi_n^*)\leq j$.
Moreover, by \eqref{2.5} and \eqref{2.11},
\begin{align*}
z\varphi_j &= \rho_j \varphi_{j+1} + \bar\alpha_j \varphi_j^* \\
\varphi_{j+1}^* &= -\alpha_j \varphi_{j+1} + \rho_j \varphi_j
\end{align*}
from which \eqref{10.4a}--\eqref{10.4d} are immediate.
\end{proof}

\begin{proof}[Third Proof of Theorem~\ref{T10.1}] This is an analog of the proof of
$\calL\calM$ factorization in Section~\ref{s2}. There $\calC$ is the matrix of overlap
of the orthonormal bases $\{z\chi_\ell\}_{\ell=0}^\infty$ and $\{\chi_\ell\}_{\ell=0}^\infty$.
The $\calL\calM$ factorization comes from inserting the basis $\{zx_\ell\}_{\ell=0}^\infty$.
Here we have the bases
\begin{align}
e^{(0)} &= (z\varphi_0, \dots, z\varphi_{N-1}) \lb{10.5} \\
e^{(N)} &= (\varphi_0, \dots, \varphi_{N-1}) \notag
\end{align}
and $\calG$ is an overlap matrix. We introduce $N-1$ intermediate bases:
\begin{align*}
e^{(1)} &= (z\varphi_0, \dots, z\varphi_{N-2}, \varphi_{N-1}^*) \\
e^{(2)} &= (z\varphi_0, \dots, z\varphi_{N-3}, \varphi_{N-2}^*, \varphi_{N-1}) \\
\cdots \\
e^{(j)} &= (z\varphi_0, \dots, z\varphi_{N-j-1}, \varphi_{N-j}^*, \varphi_{N-j+1},
\dots, \varphi_{N-1}) \\
\cdots
\end{align*}
where $e^{(N)}$ is given by \eqref{10.5} since $\varphi_0^* =1 =\varphi_0$.

Thus
\begin{equation} \lb{10.6}
\begin{aligned}
\calG_{k\ell} &= \langle e_k^{(N)},e_\ell^{(0)} \rangle \\
&= \sum_{m_1 \dots m_{N-1}}
\langle e_k^{(N)}, e_{m_1}^{(N-1)} \rangle \dots \langle e_{m_j}^{(N-j)},
e_{m_{j+1}}^{(N-j-1)}\rangle \dots \langle e_{m_{N-1}}^{(1)}, e_\ell^{(0)} \rangle
\end{aligned}
\end{equation}
is a product of $N$ matrices. $N-1$ have a change from $z\varphi_j, \varphi_{j+1}^*$ to
$\varphi_j^*, \varphi_j$ whose overlap matrix, by \eqref{10.4a}--\eqref{10.4d},
is $\wti\Theta(\alpha_j)$ and the extreme right has a change from $z\varphi_{N-1}$ to
$\varphi_{N-1}^*$, which is $\wti{\wti\Theta}(\alpha_{N-1})$ since in $L^2 (\partial\bbD)$,
\[
z\varphi_{N-1} -\bar\alpha_{N-1} \varphi_{N-1}^* =0
\]
Thus, \eqref{10.6} is \eqref{10.3}.
\end{proof}

As a first application, we want to show that each finite unitary has an $\calL\calM$
factorization without recourse to orthogonal polynomials. By taking limits, one
obtains an $\calL\calM$ factorization in general. This calculation fleshes out an
argument given by AGR \cite{Athens} in the orthogonal case by using induction
to make the proof more transparent:

\begin{theorem}\lb{T10.2A} Let $U$ be a unitary matrix on $\bbC^N$ with
$(1, 0, \dots 0)^t$ as cyclic vector. Then there exists a unitary $V$ on $\bbC^N$
with $V(1, 0, \dots 0)^t=(1, 0, \dots, 0)^t$ so that $VU\!V^{-1}$ has an $\calL\calM$
factorization.
\end{theorem}

\begin{remark} By $\calL\calM$ factorization, we mean $\calL=\wti\Theta(\alpha_0)
\oplus \wti\Theta (\alpha_2)\oplus\cdots$ and $\calM=\bdone_{1\times 1}\oplus
\wti\Theta(\alpha_1)\oplus\wti\Theta(\alpha_3)\oplus\cdots$, with a $\wti{\wti\Theta}
(\alpha_{N-1})$ at the end of $\calL$ if $N$ is odd and of $\calM$ is $N$ is even.
\end{remark}

\begin{proof} We use induction on $N$. $N=1$, which says $U=(\wti{\wti\Theta}
(\alpha_0)(\bdone)$, is trivial. By the GGT representation and AGR factorization,
we can find $W$ (with $W(1, 0, \dots, 0)^t =(1, 0,\dots, 0)^t$) so
\[
WU\!W^{-1} =\wti\Theta_0 (\alpha_0) \wti\Theta_1(\alpha_1) \dots
\wti{\wti\Theta}_{N-1}(\alpha_{N-1})
\]
Let
\begin{equation} \lb{10.6a}
U_1 = \wti\Theta_0(\alpha_1) \dots \wti{\wti\Theta}_{N-2}(\alpha_{N-1})
\end{equation}
on $\bbC^{N-1}$. By induction and adding $\bdone_{1\times 1}\oplus \cdots$
everywhere, we can find $\calL_1$, $\calM_1$, and $V_1$ so
\begin{equation}\lb{10.6b}
(1\oplus V_1)WU\!W^{-1} (1+V_1)^{-1} = \wti\Theta_0(\alpha_0)
[\bdone\oplus \calL_1][1\oplus \calM_1]
\end{equation}

Define
\[
V=[1\oplus\calM_1][1\oplus V_1]W
\]
(note $V$ maps $(1\, 0\dots 0)^t$ to itself). We have
\[
VU\!V^{-1} = \{[1\oplus\calM_1]\wti\Theta_0(\alpha_0)\}
\{[\bdone\oplus\calL_1]\}
\]
which precisely has the form $\calL\calM$.
\end{proof}

As a second application, we want to provide an explicit map that will be critical
in the next section. Fix $\delta_0\in \bbC^n$ and $U\in \bbU(n)$, the $n\times n$
unitary matrices. Let $\bbU(n-1) =\{U\in U(n)\mid U\delta_0=\delta_0\}$. The symbol
$\bbU(n-1)$ is accurate since each such $U$ defines and is defined by a unitary on
$\{\delta_0\}^\perp\cong\bbC^{n-1}$. Let $\bbS\bbC^{2n-1} = \bbU(n)/\bbU(n-1)$,
the group theoretic quotient. By mapping
\begin{equation} \lb{10.7}
\pi \colon U\in\bbU(n)\to U\delta_0
\end{equation}
we see that $\bbS\bbC^{2n-1}\cong \{z\in\bbC^n\mid\abs{z}=1\}$, the sphere of real
dimension $2n-1$. Here is the result we will need:

\begin{theorem}\lb{T10.3} There exists continuous maps $g_1$ and $g_2$ defined on
$\{z\in \bbS\bbC^{2n-1}\mid z\neq \delta_0\}$ with $g_1$ mapping to $\bbU(n)$
and $g_2$ to $\bbS\bbC^{2n-3}=\{z\in\bbS\bbC^{2n-1}\mid \langle\delta_0,
z\rangle =0\}$ so that
\begin{SL}
\item[{\rm{(i)}}]
\begin{equation} \lb{10.8}
\pi [g_1(z)] =z
\end{equation}

\item[{\rm{(ii)}}] $V(U) \equiv g_1 (\pi(U))^{-1} U\in \bbU(n-1)$ for all
$U\notin \bbU(n-1)$

\item[{\rm{(iii)}}] If $\delta_0$ is cyclic for $U$ with Verblunsky coefficients
$\alpha_j (U,\delta_0)$, then $g_2 (\pi (U))$ is cyclic for $V(U)\restriction
\bbC^{n-1}$ and
\begin{equation} \lb{10.9}
\alpha_j (V(U), g_2 (\pi(U))) = \alpha_{j+1} (U,\delta_0)
\end{equation}

\item[{\rm{(iv)}}]
\begin{equation} \lb{10.10x}
\langle\delta_0, U\delta_0\rangle = \ol{\alpha_0 (U,\delta_0)}
\end{equation}
if $\delta_0$ is cyclic for $U$\!.
\end{SL}
\end{theorem}

\begin{proof} If $z\neq \delta_0$, $a(z)=\ol{\langle \delta_0, z\rangle} \in \bbD$
and so
\[
g_2 (z) = \f{z-\langle \delta_0,z\rangle \delta_0}
{\|z-\langle \delta_0, z\rangle \delta_0\|}
\]
is well defined and in $\bbS\bbC^{2n-3}$. In particular, if $p(z)=(1-a(z)^2)^{1/2} =
\|z-\langle \delta_0, z\rangle\delta_0\|$, we have
\begin{equation} \lb{10.10}
z = p(z) g_2(z) + \ol{a(z)}\, \delta_0
\end{equation}

Define $g_1(z)$ by
\[
g_1(z) w = \begin{cases}
w \quad &\text{if } w\perp\delta_0, g_2(z) \\
z \quad &\text{if } w=\delta_0 \\
-a(z)g_2 (z) + p(z)\delta_0 \quad & \text{if } w=g_2(z)
\end{cases}
\]
and otherwise linear. Then $g_1(z)$ is unitary since $\Theta (a(z))$ is unitary.
(i) is obvious from $g_1(z)\delta_0 =z$. (ii) is a restatement of (i). (iii) follows
from the fact that $g_2(z)$ corresponds to $\delta_1$ in a $\delta_j =\chi_j (z)$
basis and the AGR factorization. (iv) is a consequence of $z\varphi_0 -\bar\alpha_0
\varphi_0^* = \varphi_1$, so $\langle\varphi_0, z\varphi_0\rangle =\bar\alpha_0
\langle\varphi_0, \varphi_0^*\rangle = \bar\alpha_0$.
\end{proof}

We want to close this section by noting that the AGR factorization implies an
estimate on the GGT matrices that is not obvious from \eqref{10.2}. Indeed, in
Section~4.1 of \cite{OPUC1}, an unnecessary condition, $\liminf \abs{\alpha_n}>0$,
is made because the estimate below is not obvious.

In essence, the AGR factorization plays the role for estimates of GGT matrices
that the $\calL\calM$ factorization does for CMV matrices. In some ways, it is
more critical because CMV matrices are five-diagonal with matrix elements which
are quadratic in $\alpha$ and $\rho$, so one can easily get estimates like
\eqref{7.2} (but with a worse constant) without using the $\calL\calM$
factorization. Since GGT matrices are not finite width and have matrix elements
that are products of arbitrary orders, direct estimates from \eqref{10.2} are
much harder.

\begin{theorem} \lb{T10.4} Let $\{\alpha_n\}_{n=0}^\infty$ and $\{\beta_n\}_{n=0}^\infty$
be two sequences in $\ol{\bbD}^\infty$ and let $\sigma_n = (1-\abs{\beta_n}^2)^{1/2}$.
Then
\begin{equation} \lb{10.11}
\|\calG (\{\alpha_n\}_{n=0}^\infty) -\calG (\{\beta_n\}_{n=0}^\infty)\|_1 \leq
2\sum_{n=0}^\infty \, (\abs{\alpha_n-\beta_n} + \abs{\sigma_n-\rho_n})
\end{equation}
\end{theorem}

\begin{proof} By \eqref{10.4} and standard trace class techniques \cite{STI}, we
need only prove for finite sequences that
\begin{equation} \lb{10.12}
\|\wti\Theta_0 (\alpha_0) \dots\wti\Theta_N(\alpha_N) - \wti\Theta_0 (\beta_0)
\dots \Theta_N (\beta_N)\|_1 \leq 2\sum_{n=0}^N \, (\abs{\alpha_n-\beta_n} +
\abs{\sigma_n -\rho_n})
\end{equation}
Since $\|\wti\Theta_j (\alpha_j)-\wti\Theta_j (\beta_j)\|_1 \leq 2(\abs{\alpha_j -\beta_j}
+ \abs{\sigma_j -\rho_j})$, and $\|\wti\Theta (\alpha)\|=1$, writing the difference of
products as a telescoping sum yields \eqref{10.12}.
\end{proof}

Notice also that AGR factorization shows that if $\abs{\alpha_j}=1$, $\calG$ decouples.
This fact and \eqref{10.11} provide an alternate proof of Rakhmanov's lemma
(Theorem~\ref{T9.2}) by the same decoupling argument, but using $\calG$ in place
of $\calC$. If $\sum_{n=0}^\infty \abs{\alpha_n}^2=\infty$, one also gets a proof
of Theorem~\ref{T9.3}, using $\calG$ in place of $\calC$ and \eqref{10.11}.
If $\sum_{n=0}^\infty \abs{\alpha_n}^2 <\infty$, one must use the extended GGT
matrix, $\calF$, of Section~4.1 of \cite{OPUC1} (see also Constantinescu
\cite{Con84}). It is easy to prove that if $\sum\abs{\alpha_n}^2 <\infty$
and $\sum\abs{\alpha_n-\beta_n}<\infty$, then $\calF (\{\alpha_n\}_{n=0}^\infty)
-\calF (\{\beta_n\}_{n=0}^\infty)$ is trace class since the difference of $\calF$'s
differs from the difference of $\calG$'s by a rank one operator, which is
always trace class!

\section{CUE, Haar Measure, and the Killip--Nenciu Theorem} \lb{s11}

In \cite{KN}, Killip and Nenciu proved the following result:

\begin{theorem}\lb{T11.1} Let $d\mu$ be a normalized Haar measure on $\bbU(n)$,
the $n\times n$ unitary matrices. Then, for a.e.\ $U$\!, $\delta_0 =(1, 0,\dots, 0)^t$
is cyclic, and the measure induced on $\bbD^{n-1} \times\partial\bbD$ by $U\to
\alpha_j (U_0,\delta_0)$ is the product measure:
\begin{equation} \lb{11.1}
\biggl\{ \prod_{j=0}^{n-2}\, \biggl[ \f{n-j-1}{\pi}\, (1-\abs{\alpha_j}^2)^{n-j-2} \,
d^2\alpha_j\biggr]\biggr\}\, \f{d\theta(\alpha_{n-1})}{2\pi}
\end{equation}
where $\theta(\alpha_{n-1})$ is defined by
\begin{equation} \lb{11.2}
\alpha_{n-1} =e^{i\theta(\alpha_{n-1})}
\end{equation}
and $d^2\alpha$ is a two-dimensional Lebesgue measure on $\bbD$.
\end{theorem}

\begin{remark} By the ``induced measure," we mean the measure $d\nu$ on $\bbD^{n-1}
\times\partial\bbD$ given by $\nu(B)=\mu(A^{-1}[B])$, where $A(U)=(\alpha_1 (U),
\dots, \alpha_{n-1}(U))$.
\end{remark}

This is really a result about Verblunsky coefficients, not CMV matrices, and
both their proof and ours use the GGT, not the CMV, matrices. We provide this
here because, first, Killip--Nenciu proved this result to provide a five-diagonal
model for CUE (see below), and because the result was proven as part of the ferment
stirred up by the CMV discovery. In this section, we will provide a partially
new proof of Theorem~\ref{T11.1} that is perhaps more natural from a group
theoretic point of view, and then describe and sketch their somewhat
shorter argument!

To understand where the factors in \eqref{11.1} come from:

\begin{lemma}\lb{L11.2A} Let $d\mu_{\bbS\bbC^{2n-1}}$ be the measure on the
$2n-1$ real dimension manifold $\{z\in\bbC^n\mid\abs{z}=1\}$, which is
normalized and invariant under rotations. Map $\bbS\bbC^{2n-1}
\overset{Q}{\longrightarrow}\ol{\bbD}$ by $z\mapsto z_1$, the first component,
and let $d\nu$ be the measure on $\ol{\bbD}$ given by $\nu(B) =
\mu_{\bbS\bbC^{2n-1}}(Q^{-1}[B])$. Then
\begin{equation} \lb{11.3}
d\nu(w) =\f{n-1}{\pi}\, (1-\abs{w}^2)^{n-2}\, d^2 w
\end{equation}
\end{lemma}

\begin{proof} Since $d^2 w=\f12 d\theta\, d\abs{w}^2$,
\[
\int_\bbD d\nu(w) = \f{2\pi(n-1)}{\pi}\, \f12 \int_0^1 (1-x)^{n-2}\, dx =1
\]
so $d\nu$ is normalized. Thus, we will not worry about constants. Using $x_1 +
ix_2, x_3 + ix_4, \dots$ for the $n$ complex variables, the measure $\delta
(1-\abs{x}^2) dx_1 \dots dx_{2n}$ is
\[
\f{dx_1\dots dx_{2n-1}}{2(1-\sum_{j=1}^{2n-1} x_j^2)^{1/2}}
\]
Integrating out $x_3, \dots, x_{2n-1}$ for fixed $x_1, x_2$ with $\rho =
(1-\abs{x_1}^2 -\abs{x_2}^2)^{1/2}$, the measure is
\[
\frac12 \int_{\abs{y}\leq (1-\rho^2)^{1/2}} \f{d^{2n-3}y}{(1-\rho^2-y^2)^{1/2}}
\]
Scaling $x=y/(1-\rho^2)^{1/2}$, we find
\[
\f12\, (1-\rho^2)^{2n-4/2} \int_{\abs{x}\leq 1} \f{d^{2n-3}x}{(1-x^2)^{1/2}}
\]
so the measure is $C(1-w^2)^{n-2} \, d^2 w$, proving \eqref{11.3}.
\end{proof}

Theorem~\ref{T11.3} below must be well known to experts on homogeneous spaces.

\begin{lemma}\lb{L11.2} Let $d\nu_1, d\nu_2$ be two probability measures on compact
spaces $X$ and $Y$, and let $d\nu$ be a probability measure on $X\times Y$\!. Suppose
\begin{SL}
\item[{\rm{(i)}}] $\pi_1^*(d\nu)=d\nu_1$, that is, if $\pi_1(x,y)=x$, then $\nu_1(B) =
\nu(\pi_1^{-1} [B])$.

\item[{\rm{(ii)}}] For any continuous $f$ on $X$\!, $\int_X f\, d\nu = C_f\, d\nu_2$,
that is,
\begin{equation} \lb{11.4}
\int f(x) g(y)\, d\nu = C_f \int g(y)\, d\nu_2
\end{equation}
for all continuous $g$ on $Y$\!.

Then $d\nu =d\nu_1\otimes d\nu_2$.
\end{SL}
\end{lemma}

\begin{proof} Taking $g=1$ in \eqref{11.4},
\[
C_f = \int f(x)\, d\nu = \int f(x)\, d\nu_1 (x)
\]
by (i). Thus,
\[
\int f(x)g(y)\, d\nu = \biggl( \int f(x)\, d\nu_1(x)\biggr)
\biggl( \int g(y)\, d\nu_2 (y)\biggr)
\]
so $d\nu= d\nu_1 \otimes d\nu_2$ integrated on product functions which are
total in $C(X\times Y)$.
\end{proof}

\begin{theorem}\lb{T11.3} Let $G$ be a compact group and $H$ a closed subgroup.
Let $d\nu_G, d\nu_H$ be normalized Haar measures and $\pi\colon G\to G/H$.
Let $d\nu_{G/H}$ be the measure induced by $d\nu_G$ on $G/H$, that is,
\[
\nu_{G/H}(B) = \nu_G(\pi^{-1}[B])
\]
Let $\calO$ be an open set in $G/H$ and $f\colon\calO\to G$ a continuous
cross-section, that is, $\pi [f(x)]=x$ for all $x\in\calO$. Coordinatize
$\pi^{-1}[\calO]$ by $\calO\times H$ via
\begin{equation} \lb{11.5}
(x,h)\mapsto f(x) h
\end{equation}
Then, on $\pi^{-1} [\calO]$,
\begin{equation} \lb{11.6}
d\nu_G(x,h) = d\nu_{G/H}(x)\, d\nu_H (h)
\end{equation}
\end{theorem}

\begin{proof} The existence of a cross-section implies that under the coordinates
\eqref{11.5}, $\pi^{-1}[\calO]\cong\calO\times H$. Clearly, $\pi_1^* (d\mu_G)=
d\mu_{G/H}\restriction \calO$, by construction of $d\mu_{G/H}$. On the other
hand, $\int_\calO f\, d\mu_{G/H}$ is a measure on $H$ invariant under right
multiplication by any $h\in H$\!, so this is $C_f\, d\mu_H$. Therefore,
Lemma\ref{L11.2} applies and \eqref{11.6} holds.
\end{proof}

\begin{proof}[Proof of Theorem~\ref{T11.1}] We use induction on $n$. $n=1$,
that is, that for $\bbU(1)$, $U=(e^{i\theta_0})$ has Haar measure
$\f{d\theta_0}{2\pi}$, is immediate.

Note that $U\in\bbU(n)$ has $\delta_0$ as a cyclic vector if and only if $U$ has
simple spectrum, and for each eigenvector $\varphi_k$ of $U$\!, we have
$\langle\varphi_k,\delta_0\rangle \neq 0$. As is well known, $\{U\mid U$
has a degenerate eigenvalue$\}$ has codimension $3$ and so zero Haar measure.
Similarly, $\langle\varphi_k,\delta_0\rangle =0$ on a set of codimension $2$
and so zero Haar measure. Thus, $\calC_n =\{E\mid\delta_0$ is cyclic for
$U\}$ has full Haar measure.

Let $\calO=\{\eta\in \bbS\bbC^{2n-1}\mid\eta\neq\delta_0\}$. Then $f(x)=g_1(z)$
given in Theorem~\ref{T10.3} is a cross-section, and so $d\mu_{\bbU(n)}=
d\mu_{\bbS\bbC^{2n-1}}\otimes d\mu_{\bbU(n-1)}$ by Theorem~\ref{T11.3}.

By Theorem~\ref{T10.3}, $g_1^{-1} (\pi_1 [\calC_n])$ is $(z,V(U))$ and
$V(U)$ has Verblunsky coefficients $\{\alpha_{j+1}(U)\}_{j=0}^{n-2}$.
Thus, by induction, $d\mu_{\bbU(n-1)}$ on these $\alpha$'s is the product
\eqref{11.6} without the $\alpha_0$ factor.

By \eqref{10.10} and Lemma~\ref{L11.2}, the $\alpha_0$ distribution generated
by $d\mu_{\bbS\bbC^{2n-1}}$ is the $\alpha_0$ factor in \eqref{11.1}.
\end{proof}

The proof in \cite{KN} differs in two ways: First, in place of the AGR factorization,
Killip--Nenciu use a (Householder) factorization as a phase factor times a product of
reflections. Instead of using induction on symmetric spaces as we do, they use
an alternate that would work with the AGR factorization also. Starting with
$\varphi_0=\delta_0$, we let $\psi_0=U\varphi$. There is a unique vector,
$\varphi_1$ (what we called $g_2(\pi(U))$ in Theorem~\ref{T10.3}), in the span
of $\varphi$ and $\psi_0$, so that $\langle\psi_0,\varphi_1\rangle >0$ and
$\langle \varphi_1, \delta_0\rangle =0$. $\varphi_1$ is cyclic for $V(U)$ and so,
by induction, we obtain $\psi_0,\psi_1,\dots,\psi_{n-1}$ an $O\!N$ basis with
$\abs{\langle\delta_j,\psi_j\rangle} <1$. It is not hard to see that, via the
AGR factorization, this sets up a one-one map of $O\!N$ basis with $\abs{\langle
\delta_j, \psi_j\rangle} <1$ and $U$'s with $\delta_0$ cyclic for $U$\!. Haar
measure induces on the $\psi$'s a measure as follows: $\psi_0$ is uniformly
distributed on $\bbS\bbC^{2n-1}$; $\psi_1$ uniformly on the copy of $\bbS\bbC^{2n-3}$
of unit vectors orthogonal to $\psi_1$; $\psi_2$ uniformly on $\bbS\bbC^{2n-5}$,
etc. Since $\langle\delta_j,\psi_j\rangle =\bar\alpha_j$, we obtain the measure
\eqref{11.1}.

Since $\det (\wti\Theta_j(\alpha_j))=-1$, $\det \bigl(\wti{\wti\Theta}_{N-1}
(\alpha_{N-1})\bigr)=\bar\alpha_{N-1}$, we see
\[
\det (\calG (\{\alpha_n\}_{n=0}^{N-1}) = (-1)^{N-1} \bar\alpha_{N-1}
\]
Thus, $\bbS\bbU(n) =\{U\in\bbU(n)\mid\det(U)\}$ is precisely these $U$ with
$\alpha_{N-1}=(-1)^{N-1}$. The same inductive argument thus proves:

\begin{theorem}\lb{T11.4} Let $d\mu$ be normalized Haar measure on
$\bbS\bbU(n)$. Then for a.e.\ $U$\!, $\delta_0=(1, 0,\dots, 0)^t$ is cyclic
and the measure induced on $\bbD^{n-1}$ by $U\to \alpha_j (U,\delta_0)$
{\rm{(}}with $\alpha_{n-1}(U,\delta_0)\equiv (-1)^{n-1}${\rm{)}} is the
product measure given by \eqref{11.1} with the final $d\theta$ term dropped.
\end{theorem}

$\bbS\bbO[n]$ is the $n\times n$ real unitary matrices (i.e., orthogonal
matrices). If $\delta_0$ is cyclic, they have Verblunsky coefficients which are
easily seen to lie in $(-1,1)$. Conversely, it is easy to see that if $\alpha_j\in
(-1,1)$ for $j=0,\dots, n-2$, there is an orthogonal matrix with those $\alpha_j$'s.
A similar analysis lets us compute the distribution on $(-1,1)^{n-1}$ induced
by Haar measure on $\bbS\bbO[n]$. We need only replace Lemma~\ref{L11.2} by

\begin{lemma}\lb{L11.5} Let $d\eta_{n-1}$ be the measure on the $n-1$-dimensional
unit sphere in $\bbR^n$. The induced measure on $x_1$ is
\[
\f{\Gamma(\f{n}{2}) (1-\abs{x_1}^2)^{(n-3)/2}\, dx_1}{\sqrt{\pi}\, \Gamma (\f{n-1}{2})}
\]
\end{lemma}

\begin{proof} That the measure is $C(1-\abs{x_1}^2)^{(n-3/2)}\, dx_1$ follows from
the same calculation as in Lemma~\ref{L11.2}. The normalization is the inverse
of the beta function $2^{2-n}\Gamma(n-1)/\Gamma (\f{n-1}{2})$ which, as noted
by \cite{KN} (there is a $(\,\,)^{-1}$ missing on the leftmost term in
their (3.14)), can be written, using the duplication formula for beta functions,
as $\Gamma (\f{n}{2})/\sqrt\pi\, \Gamma (\f{n-1}{2})$.
\end{proof}

We thus have

\begin{theorem}[\cite{KN}] \lb{T11.6} The measure on $(-1,1)^{n-1}$ induced
by Haar measure on $\bbS\bbO[n]$ mapped to the real Verblunsky coefficients is
\[
\f{\Gamma (\f{n}{2})}{\pi^{n/2}}\, \prod_{k=0}^{n-1}
(1-\alpha_k^2)^{(n-k-3)/2}\, d\alpha_k
\]
\end{theorem}

The CUE eigenvalue distribution \cite{Dyson1,Dyson2,Dyson3} is the one for
$U\in \bbU[n]$ induced by Haar measure. Weyl's integration formula (see, e.g.,
\cite{SGR}) says that if $\lambda_1, \dots, \lambda_n$ with $\lambda_j =
e^{i\theta_j}$ are the eigenvalues, this is $C \prod_{i<j} \abs{\lambda_i -
\lambda_j}^2 \prod_{j=1}^n \f{d\theta_j}{2\pi}$. Theorem~\ref{T11.1} says
that CMV matrices with distribution of $\alpha$'s given by \eqref{11.1}
has the same distribution, and so gives a model for CUE by five-diagonal
matrices. \cite{KN} find a similar model for the ``$\beta$-distributions"
given by $C_\beta \prod_{i=j} \abs{\lambda_i - \lambda_j}^\beta
\f{d\theta_j}{2\pi}$; see also Forrester--Rains \cite{FR}.

\section{CMV and the AL Flow} \lb{s12}

One of the great discoveries of the 1970's (\cite{DubMatNov,FlMcL,GGKM,McvM1,vMoer}
and dozens of other papers) is that lurking within one-dimensional Schr\"odinger
operators and Jacobi matrices is a completely integrable system (resp., KdV and
Toda flows), natural ``invariant" tori, and a natural symplectic structure in
which the Schr\"odinger operator or Jacobi matrix is the dynamical half of a Lax
pair.

Such structures occur also for Verblunsky coefficients, and the dynamical half
of the Lax pair is the CMV matrix. While the CMV part obviously requires CMV
matrices, the other parts do not, and it is surprising that it was only in
2003--04 that they were found. We will settle here for describing the two
most basic structures, leaving further results to mentioning the followup
papers: Geronimo--Gesztesy--Holden \cite{GGH}, Gesztesy--Zinchenko \cite{GZ},
Killip--Nenciu \cite{KNnew}, Li \cite{Li}, and Nenciu \cite{NenRev}.

On $\bbD$, introduce the symplectic form given by the Poisson bracket
(where, as usual, $\f{\partial}{\partial z}$ and $\f{\partial}{\partial \bar z}$
are $\f12[\f{\partial}{\partial x} \mp i\f{\partial}{\partial y}]$),
\begin{equation} \lb{12.1}
\{f,g\} = i\rho^2 \biggl[ \f{\partial f}{\partial\bar z} \, \f{\partial g}{\partial z}
- \f{\partial f}{\partial z} \, \f{\partial g}{\partial \bar z}\biggr]
\end{equation}
The $\rho^2$ is natural as we will see below. Extend this to $\bbD^p$
(coordinatized by $(\alpha_0, \dots, \alpha_{p-1})$ by
\begin{equation} \lb{12.2}
\{f,g\} = i \sum_{j=0}^{p-1} \rho_j^2 \biggl[ \f{\partial f}{\partial\bar\alpha_j} \,
\f{\partial g}{\partial\alpha_j } - \f{\partial f}{\partial\alpha_j} \,
\f{\partial g}{\partial \bar\alpha_j}\biggr]
\end{equation}
Because of the $\rho^2$,
\[
\biggl\{\,\prod_{j=0}^{p-1} \rho_j^2, g\biggr\} =
-\prod_{j=0}^{p-1} \rho_j^2 \biggl(\, \sum_{j=0}^{p+1}
\f{\partial g}{\partial\theta_j}\biggr)
\]
and functions of $\rho_0 \dots \rho_{p-1}$ generate simultaneous rotations of all
phases. Nenciu--Simon \cite{NS} proved the following:

\begin{theorem}[\cite{NS}]\lb{T12.1} Let $p$ be even and let $\Delta(z,
\{\alpha_j\}_{j=0}^{p-1})$ be the discriminants {\rm{(}}see \eqref{3.15}{\rm{)}}
for the periodic sequence with $\alpha_{j+kp}=\alpha_j$ for $j=0,\dots, p-1$;
$k=0,1,2,\dots$. Then, with respect to the symplectic form \eqref{12.2},
\begin{equation} \lb{12.3}
\{\Delta(w), \Delta(z)\}=0
\end{equation}
for all $w,z\in\bbC\backslash\{0\}$.
\end{theorem}

{\it Note}: See \cite{Li,KNnew} for a discussion of symplectic forms on unitary
matrices.

\smallskip
Since $\ol{\Delta (1/\bar w)}=\Delta (w)$, and the leading coefficient is real,
$\Delta(z)$ has $p$ real coefficients, that is, $\Delta (z) =\sum_{j= -p/2}^{p/2}
a_j z^j$ with $a_{-j}=\bar a_j$, then $a_{p/2}, \Real a_{p/2-1}, \Ima a_{p/2 -1},
\dots, \Real a_1, \Ima a_1, a_0$ are the $p$ real functions of $\{\alpha_j\}_{j=0}^{p-1}$
which Poisson commute. They are independent at a.e.\ points (in $\alpha$) and
define invariant tori. Each one generates flows that are completely integrable.
The simplest is
\begin{equation} \lb{12.4}
-i\dot{\alpha}_j =\rho_j^2 (\alpha_{j+1} + \alpha_{j-1})
\end{equation}
which has been known as a completely integrable system for a long time under
the name ``defocusing Ablowitz--Ladik flow" (after \cite{AL75,AL76b,AL76a}).

Nenciu has proven a beautiful result:

\begin{theorem}[\cite{Nendiss,NenIMRN,NenRev}] The flows generated by the coefficients
of $\Delta$ can be put into Lax pair form with the dynamical element being the Floquet
CMV matrix.
\end{theorem}

For details as well as extensions to some infinite CMV matrices, see the
references above.

The flow generated by $\prod_{j=0}^{p-1} \rho_j^2$ realizes the $\alpha_j \to
\lambda\alpha_j$ invariance of the isospectral tori. The flow \eqref{12.4} is
generated by $\Real (a_1)$. The $\Ima (a_k)$ generate flows that preserve the
set of $\{\alpha_j\}_{j=0}^{p-1}$ where all $\alpha_j$ are real (as a set,
not pointwise). The simplest of these, generated by $\Ima (a_1)$, is
\begin{equation} \lb{12.5}
\dot{\alpha}_n = \rho_n^2 (\alpha_{n+1} -\alpha_{n-1})
\end{equation}
called the Schur flow. Via the Geronimus relations of the next section, these generate
a flow on Jacobi parameters that is essentially the Toda flow. For further discussion,
see \cite{AG94,FG,Gek,Nendiss}.

\section{CMV Matrices and the Geronimus Relations} \lb{s13}

In a celebrated paper, Szeg\H{o} \cite{Sz22a} found a connection between orthogonal polynomials
for measures on $[-2,2]$ (he had $[-1,1]$; I use the scaling common in the Schr\"odinger
operator community) and OPUC. Given a measure $d\gamma$ on $[-2,2]$, one defines the
unique measure $d\xi$ on $\partial\bbD$ which is invariant under $z\to\bar z$ and
obeys
\begin{equation} \lb{13.1}
\int g(x)\, d\gamma(x) = \int g(2\cos\theta)\, d\xi (\theta)
\end{equation}
What Szeg\H{o} showed is that the orthonormal polynomials $p_n$ for $d\gamma$ and the
OPUC for $\varphi_n$ for $d\xi$ are related by
\begin{equation} \lb{13.2}
p_n \biggl( z+\f{1}{z}\biggr) = C_n z^{-n} (\varphi_{2n}(z) + \varphi_{2n}^*(z))
\end{equation}
The normalization constants $C_n$ (see (13.1.14) in \cite{OPUC2}) $\to 1$ as $n\to\infty$
if $\sum_{n=0}^\infty \abs{\alpha_n}^2 <\infty$. Motivated by this, Geronimus \cite{Ger46}
found a relation between the Verblunsky coefficients, $\alpha_n$, for $d\xi$ and the
Jacobi parameters $\{a_n\}_{n=1}^\infty, \{b_n\}_{n=1}^\infty$ for $d\gamma$ (see
Theorem~13.1.7 of \cite{OPUC2}):

\begin{subequations} \lb{13.3}
\begin{align}
a_{n+1}^2 &= (1-\alpha_{2n-1}) (1-\alpha_{2n}^2)(1+\alpha_{2n+1}) \lb{13.3a} \\
b_{n+1} &= (1-\alpha_{2n-1}) \alpha_{2n} - (1+\alpha_{2n-1})\alpha_{2n-2} \lb{13.3b}
\end{align}
\end{subequations}

In \cite{KN}, Killip--Nenciu found a direct proof of \eqref{13.3} by finding a
beautiful relation between CMV and some Jacobi matrices. We will sketch the idea,
leaving the detailed calculations to \cite{KN} or the pedagogic presentation
in Section~13.2 of \cite{OPUC2}.

A measure is invariant under $z\to\bar z$ if and only if all $\{\alpha_n\}_{n=0}^\infty$
are real. $\Theta(\alpha)$ with $\alpha$ real is selfadjoint and unitary with determinant
$-1$, hence eigenvalues $\pm 1$, that is, a reflection on $\bbC^2$. Thus,
\[
\alpha_n =\bar\alpha_n \, \text{ all } n \Rightarrow \calM^2 = \calL^2 =\bdone
\]
Since $\chi_n(z) = \ol{x_n (1/\bar z)}$, we see that if $\mu$ is invariant and
$(Mf)(z) = f(\bar z)$, then
\[
\langle \chi_j, M\chi_\ell\rangle = \calM_{j\ell}
\]

$\calC + \calC^*$ is selfadjoint and maps $\{f\mid \calM f = f\}$ to itself.
Let us see in a natural basis that its restriction to this invariant subspace
is a Jacobi matrix.

If $\alpha$ is real and
\[
S(\alpha) =\f{1}{\sqrt 2}\, \begin{pmatrix}
\sqrt{1-\alpha} & -\sqrt{1+\alpha} \\
\sqrt{1+\alpha} & \sqrt{1-\alpha} \end{pmatrix}
\]
then
\begin{equation} \lb{13.4}
S(\alpha)\Theta (\alpha) S(\alpha)^{-1} =
\begin{pmatrix}
-1 & 0 \\
0 & 1
\end{pmatrix}
\end{equation}
Define
\[
\calS = \bdone_{1\times 1} \oplus S(\alpha_1)\oplus S(\alpha_3 ) \oplus \cdots
\]
so
\[
\calS\calM \calS^{-1} =\calR =
\begin{pmatrix}
{} & 1 \\
{} & {} & -1  \\
{} & {} & {} & 1   \\
{} & {} & {} & {} & -1 &  \\
{} & {} & {} & {} & {} & \ddots & {}
\end{pmatrix}
\]
and define
\[
\calB = \calS\calL\calS^{-1}
\]
Then
\[
\calS(\calC + \calC^{-1})\calS^{-1} = \calR\calB + \calB\calR
\]
which commutes with $\calR$.

$\calB$ is seven-diagonal as a product of three tridiagonal matrices. Moreover,
since $\calB$ commutes with $\calR$, its odd-even matrix elements vanish. It
follows that
\[
\calR\calB + \calB\calR = \calJ_e \oplus\calJ_o
\]
where $\calJ_e$ acts on $\{\delta_{2n}\}_{n=0}^\infty$ and $\calJ_o$ on
$\{\delta_{2n+1}\}_{n=0}^\infty$, and each is a Jacobi matrix. A calculation
shows that the Jacobi parameters of $\calJ_e$ are given by \eqref{13.3}
and that the spectral measures are related by \eqref{13.1}. One can also
analyze $\calJ_o$ which is related to another mapping of Szeg\H{o}
\cite{Sz22a} and one gets two more Jacobi matrices by looking at $\calC +
\calC^{-1}$ restricted to the spaces where $\calL =1$ or $\calL=-1$.

\bigskip

\end{document}